%% file: gen_lqr.tex
\documentclass[12pt]{article}
\usepackage{fullpage}
\usepackage{graphicx,verbatim,booktabs,amsmath,amsfonts,colortbl,xcolor,amsthm,siunitx}
\usepackage{hyperref}
\usepackage[small,bf]{caption}
\usepackage{subcaption} 

\input{defs.tex}

\begin{document}

\date{}

\title{Stochastic Control with Affine Dynamics and Extended Quadratic Costs}

\author{
Shane Barratt \\
Stanford University \\
sbarratt@stanford.edu
\and
Stephen Boyd \\
Stanford University \\
boyd@stanford.edu
}

\maketitle

\begin{abstract}
An extended quadratic function is a quadratic function plus the indicator function
of an affine set, \ie, a quadratic function with embedded linear equality constraints.
We show that, under some technical conditions, 
random convex extended quadratic functions are 
closed under addition, composition with an affine function, expectation, and 
partial minimization, \ie, minimizing over some of its arguments.
These properties imply that dynamic programming can be tractably 
carried out for stochastic control problems with random affine dynamics and 
extended quadratic cost functions.
While the equations for the
dynamic programming iterations are
much more complicated than for traditional linear quadratic control,
they are well suited to an object-oriented implementation, which we describe.
We also describe a number of known and new applications.
\end{abstract}

\newpage
\tableofcontents
\newpage

\section{Introduction}
\label{sec:introduction}

Many practical problems can be modeled as stochastic control problems.
Dynamic programming, pioneered by Bellman in the 1950's~\cite{bellman1954theory}, 
provides a solution method, at least in principle~\cite{bertsekas2017dynamic}.
Dynamic programming relies on the cost-to-go, value, or Bellman function (on the state space),
which is computed by an iteration involving a few operations such as addition,
expectation over random variables, and minimization over the allowed actions
or controls.
The cost-to-go function can be tractably represented, and 
these operations can be carried out tractably in only a few special settings.
\begin{itemize}
\item \emph{Finite state and control spaces.}
In this case, the cost-to-go function and the control policy 
can be explicitly represented by lookup tables.
\item \emph{Vector-valued states and controls, linear dynamics, and convex quadratic cost.}
This is the famous linear quadratic control or regulator (LQR) problem.
In this case, the cost-to-go function is a convex quadratic form, represented
by a matrix, and the control policy is linear, represented by a gain 
matrix~\cite{kalman1964linear}.
\end{itemize}
Despite the special forms of these two cases, they are very widely applied.
There are a few other very specialized cases where 
dynamic programming is tractable, such as the optimal consumption 
problem~\cite{merton1971optimum}.

For cases where exact dynamic programming is intractable, 
many methods have been developed to approximately solve the problem, such as
approximate dynamic programming (ADP), reinforcement learning, and many others.
These methods can be very effective in practice, depending on the approximations
or algorithms used, which can vary across applications.
There is a vast literature on these methods;
see, \eg,~\cite{bertsekas1988adaptive,bertsekas1996neuro,sutton1998reinforcement,powell2007approximate,wang2009performance,wang2011performance,keshavarz2012quadratic,o2013iterated},
and the many references in them.

Our focus in this paper is to
identify a class of stochastic control problems that, like the two special cases above,
can be solved exactly. Our class is a generalization of the classic LQR problem.
The class of problems, which we formally describe in the next section, 
has a state with a vector-valued and finite part (which we call the mode),
a vector-valued control,  
random mode-dependent affine dynamics, random mode-dependent extended 
quadratic cost, and state/control-independent Markov chain dynamics for the mode.
Extended quadratic functions, which we define formally in the section section,
are quadratic functions that include linear and constant terms, as well as 
implicit linear equality constraints.
Many special cases of our general problem class have been noted and solved in the 
literature, \eg, so-called jump-linear quadratic control \cite{florentin1961optimal,krasovsky1961analytical} and LQR with random dynamics~\cite{drenick1964optimal}.
We unify these problems under one common problem description and solution method;
in addition, our class includes problems that, to the best of our knowledge, 
have not been addressed in the literature.

While dynamic programming for our class of problems can be carried out exactly
(modulo how expectation is carried out), the equations that characterize the 
cost-to-go function and the policies are not simple, and in particular, are far 
more complex than those for LQR, which are well known.
Our approach is to develop an object-oriented solution method.
To do this, we identify the key functions and methods
that must be carried out, and describe how to implement them; dynamic programming
simply uses these methods, without expanding the equations and formulas.
This approach has several advantages. First, it can be immediately implemented
(and indeed, has been).
Second, it focuses on the critical ideas without getting bogged down 
in complicated equations, as the traditional approach would.
Third, its generality and compositional form allows it to apply to a wide variety of 
problems; in particular, components can be re-arranged to solve other 
problems not described here.

\subsection{Related work}

Stochastic control has applications in a wide variety of areas, including supply-chain optimization~\cite{arrow1951optimal,sarimveis2008dynamic}, advertising~\cite{nerlove1962optimal}, finance~\cite{merton1971optimum,samuelson1975lifetime}, dynamic resource allocation~\cite{elmaghraby1993resource,bertsekas1999rollout}, and traditional 
automatic control \cite{anderson2007optimal}.
Dynamic programming is by far the most commonly employed solution method for stochastic control problems.
Dynamic programming was pioneered by Bellman in the 1950's~\cite{bellman1954theory}; for a modern treatment and its applications to stochastic control see the textbooks by Bertsekas~\cite{bertsekas1978stochastic,bertsekas2017dynamic,bertsekas2012dynamic} and the many references in them.

The linear quadratic Gaussian (LQG) stochastic control problem traces back to Kalman in the late 1950's~\cite{kalman1960contributions} and was studied
heavily throughout the 1960's (see~\cite{ieee1971} and its many references).
Since then, many tractable extensions of LQG have been proposed.
Two of the most notable extensions that have been formulated (and solved) are jump linear quadratic control (jump LQR) and random LQR, which we now describe.

\paragraph{Jump LQR.}
One special case of the class of problems described in this paper is jump LQR, where the dynamics are linear but suddenly change according to a fully observable Markov chain process.
The optimal policy for this problem in continuous time was first identified by Krasovsky and Lidsky in 1961~\cite{krasovsky1961analytical} and Florentin~\cite{florentin1961optimal}, and discovered independently by Sworder in 1969~\cite{sworder1969feedback}.
The problem was then solved in discrete time~\cite{blair1975feedback}, where the authors found that the cost-to-go functions are quadratic for each mode and that the optimal policy is linear for each mode.
These results were extended to the infinite-horizon case by Chizeck~\etal~\cite{chizeck1983markovian} and to have equality constraints in the cost by Costa~\etal~\cite{costa1999constrained}.
(See~\cite{costa2006discrete} and the references therein for a comprehensive overview of jump LQR.)
Jump LQR was applied early on to robust control system design~\cite{birdwell1977reliability}, later to reliable placement of control systems components~\cite{montgomery1985reliability}, and also to distributed control with random delays~\cite{nilsson1998real}.

\paragraph{Random LQR.}
The other important special case is random LQR.
In random LQR, the goal is to control a system that has random affine dynamics and quadratic stage cost.
This problem was first identified and solved by Drenick and Shaw in 1964~\cite{drenick1964optimal} and then in continuous time by Wonham in 1970~\cite{wonham1970random}.
For a more modern treatment in discrete time, see the paragraph ``Random System Matrices'' on pages 123-124 in~\cite{bertsekas2017dynamic}.
The random LQR problem was extended to have (jointly) random quadratic stage cost and equality constraints in~\cite{boyd2012stochastic}, and has been applied to finance in~
\cite{bertsimas1998optimal,bertsimas1999optimal} and~\cite{busetti2018thesis}.
All of these works have (somewhat independently) derived that the cost-to-go functions are quadratic and that the optimal policies are an affine function of the state.

To the best of our knowledge, no one has combined the jump and random LQR problems into a general problem class and identified the form of the solution to this general problem class.
This paper can be viewed as a unification of these two problem classes, while still maintaining the familiar tractability of LQR.

\subsection{Problem statement}

In this section we formally describe the class of stochastic control problems that we consider.

\paragraph{Random jump linear dynamical systems.}

We consider discrete-time dynamical systems, with dynamics
described by
\BEQ
\begin{split}
x_{t+1}&= f_t^{s_t}(x_t,u_t,w_t),
\quad t=0,1,\ldots,\\
s_{t+1}&=\text{$i$ with probability $\Pi_{t,ij}$ if $s_t=j$}, \quad t=0,1,\ldots,
\end{split}
\label{eq:dynamics}
\EEQ
where $t$ indexes time.
Here $x_t \in\reals^n$ (the set of real $n$-vectors) is the \emph{state} of the system at time $t$,
$s_t \in \{1,\ldots,K\}$ is the \emph{mode} of the system at time $t$,
$u_t \in \reals^m$ is the \emph{control} or \emph{input} to the system at time $t$,
$w_t \in \mathcal{W}_t$ is a random variable corresponding to the \emph{disturbance} at time $t$,
$f_t^s:\reals^n \times \reals^m \times \mathcal{W}_t \rightarrow
\reals^n$ are the \emph{state transition functions} at time $t$ when the system is in mode $s$,
and $\Pi_t$ is the \emph{mode switching probability matrix} at time $t$.

In this paper we consider state transition functions that are affine
in $x$ and $u$, \ie,
\[
f_t^s(x,u,w) = A_t^s(w)x + B_t^s(w)u + c_t^s(w), \quad t=0,1,\ldots,
\]
where 
$A_t^s: \mathcal{W}_t \rightarrow \reals^{n\times n}$ (the set of real $n\times n$ matrices)
is the \emph{dynamics matrix} at time $t$ when the system is in mode $s$,
$B_t^s: \mathcal{W}_t \rightarrow \reals^{n\times m}$
is the \emph{input matrix} at time $t$ when the system is in mode $s$,
and
$c_t^s: \mathcal{W}_t \rightarrow \reals^{n}$
is the \emph{offset} at time $t$ when the system is in mode $s$.

\paragraph{Independence assumptions.}
Because the disturbances $w_t$ are random variables, this makes $x_t$, $u_t$, and $s_t$ all random variables.
We assume that $w_t$ is independent of $x_t$, $u_t$, $s_t$, and 
$w_{t'}$ for $t' \neq t$.
We often have that $f_t$, $\Pi_t$, and the distribution of $w_t$ do not depend on $t$,
in which case the dynamics are said to be \emph{time-invariant}.
In some applications, the dynamics matrix, the input matrix, or the offset do not depend
on $w_t$, \ie, they are deterministic.

\paragraph{Information pattern.}

At time $t$, we choose $u_t$ given knowledge of the previous states $x_0,\ldots,x_t$ and modes $s_0,\ldots,s_t$, but no knowledge of the disturbance $w_t$.
For the problem we consider it can be shown that there is an optimal policy that only depends on the current state and mode~\cite{bertsekas2017dynamic}, \ie, we can express an optimal policy as
\[u_t = \phi_t^{s_t}(x_t), \quad t=0,1,\ldots,\]
where $\phi_t^s: \reals^n \rightarrow \reals^m$ is called the
\emph{policy} at time $t$ for mode $s$.
When we refer to $\phi_t$ without the superscript, we are referring to
the collection of policies at that time step.
If $\phi_t$ do not depend on $t$, then the policy is
said to be \emph{time-invariant} and is denoted $\phi$.

\paragraph{Finite-horizon problem.}

In the \emph{finite-horizon problem}, our objective is to
find a sequence of policies $\phi_0,\ldots,\phi_{T-1}$ that minimize the
expected cost over a finite time horizon, given by
\begin{equation}
J(\phi_0,\ldots,\phi_{T-1})=
\Expect
\left[\sum_{t=0}^{T-1} g_t^{s_t}(x_t,\phi_t^{s_t}(x_t),w_t) + g_T^{s_T}(x_T)\right],
\label{eq:cost-finite}
\end{equation}
subject to the dynamics~\eqref{eq:dynamics}, where $T$ is the horizon
length, the expectation is over $w_0,\ldots,w_{T-1}$ and $g_t^s:
\reals^n\times\reals^m\times\mathcal{W}_t
\rightarrow
\reals\cup\{+\infty\}$
is the \emph{stage cost function} for time $t$ when the system is in mode $s$,
and $g_T^s: \reals^n \rightarrow \reals$
is the \emph{final} stage cost function.

In this paper we consider stage cost functions of the form
\[
g_t^s(x,u,w)=\frac{1}{2}\begin{bmatrix}x \\ u \\ 1\end{bmatrix}^T G_t^s(w)
\begin{bmatrix}x \\ u \\ 1\end{bmatrix} +
\begin{cases}
    0 & F_t^s x + H_t^s u + h_t^s = 0 \\
    +\infty & \text{otherwise}
\end{cases},
\]
where $G_t^s: \mathcal{W}_t \rightarrow \symm^{n+m+1}$ ($\symm^n$ denotes the set of real symmetric $n \times n$ matrices),
$F_t^s \in \reals^{p \times n}$, $H_t^s \in \reals^{p \times m}$, and
$h_t^s \in \reals^p$.
Stage cost functions that have this form are 
\emph{extended quadratic} functions of $x$ and $u$, 
which are quadratic functions with embedded linear equality constraints,
discussed in much greater detail in \S\ref{sec:preliminaries}.
We consider a final stage cost function of the form
\[
g_T^s(x)=
\frac{1}{2}
\begin{bmatrix}
    x \\ 1
\end{bmatrix}^T
G_T^s
\begin{bmatrix}
    x \\ 1
\end{bmatrix}
+
\begin{cases}
    0 & F_T^s x + h_T^s = 0 \\
    +\infty & \text{otherwise}
\end{cases},
\]
where $G_T^s \in \symm^{n+1}$, which is an extended quadratic function of $x$.

\paragraph{Infinite-horizon.}

In the \emph{infinite-horizon problem}, we assume that the dynamics and
cost are time-invariant.
Our goal is to find a time-invariant policy
$\phi$ that minimizes the expected cost over an infinite time horizon, given by
\begin{equation}
J(\phi)=
\lim_{T\rightarrow\infty}
\Expect
\left[
\sum_{t=0}^T \gamma^t g^{s_t}(x_t,\phi^{s_t}(x_t),w_t)
\right],
\label{eq:cost-infinite}
\end{equation}
subject to the dynamics \eqref{eq:dynamics}, where $\gamma \in (0,1]$
is the discount factor, the expectation is taken over $w_t$,
and $g^s:\reals^n \times \reals^m \times \mathcal{W}_t\rightarrow \reals \cup \{+\infty\}$ is the stage cost function indexed by $s$.
We consider stage cost functions that are extended quadratic functions of $x$ and $u$.
One can recover the \emph{undiscounted} infinite-horizon problem by letting
$\gamma=1$.

We call the above problems \emph{extended quadratic control} problems.

\paragraph{Problem data representation.}

Throughout the paper, we assume that we have, at a bare minimum, access to an \emph{oracle} that provides independent samples of the random quantities
\[
A_t^s(w),\; B_t^s(w),\; c_t^s(w),\; G_t^s(w),
\]
for all $t,s$.
The samples can be given in batch, \eg, a sample of $N$ dynamics
matrices $A_t^s$ for time $t$ would result in a $N \times K \times n \times n$
matrix.
We will see later that, in some cases (namely, when the cost-to-go function is a non-extended quadratic function), additional knowledge of the distributions (in particular, their first and second moments) can be used to derive analytic expressions for expectations of quadratic functions.

In addition, we assume that we have access to the (deterministic) quantities
\[
\Pi_t, \; F_t^s, \; H_t^s, \; h_t^s, \; G_T^s, \; F_T^s, \; h_T^s,
\]
for all $t,s$.
These could be represented by matrices,
\eg, $\Pi_t$ for time $t$ would be a
$K \times K$
matrix.

\paragraph{Pathologies.}
There are several pathologies that can (and often do) occur in our formulation,
depending on the exact
problem data and distributions.
\begin{itemize}
	\item \emph{Infinite cost.} This happens if, \eg, the equality constraints are impossible to satisfy or the expectation is $+\infty$
for all policies.
	\item \emph{Cost that is unbounded below.} There exist policies that achieve arbitrarily low cost.
	\item \emph{Cost that is undefined.} 
The expectations in \eqref{eq:cost-finite} or \eqref{eq:cost-infinite} do 
not exist.
\end{itemize}
Many of these pathologies are discussed in great detail in~\cite{bertsekas2012dynamic} and~\cite{bertsekas1978stochastic}.

In this paper we do not focus on analyzing when these pathologies occur in the class of problems that we consider,
but rather on the practical application and implementation of these methods.
Also, in well posed practical problems, these pathologies rarely occur.
Nevertheless, the algorithms that we describe are capable of catching many of
these pathologies and reporting the nature of the pathology.

\subsection{Results}

In the absence of the pathologies described above, we show in this paper that there is an optimal policy in the finite-horizon problem that is an affine function of $x$, meaning the policy has the form
\BEQ
\phi_t^s(x)=K_t^sx + k_t^s,
\label{eq:policies}
\EEQ
where $K_t^s \in \reals^{m\times n}$ is the \emph{input gain matrix} and
$k_t^s \in \reals^{m}$ is the \emph{input offset matrix}.
For the infinite-horizon problem, there is a policy that is a time-invariant affine function of $x$. 
Also, the cost-to-go functions are extended quadratic functions of $x$ for each mode $s$.

When $k_t^s \in \range(K_t^s)$ (the range of the matrix $K_t^s$), we can express the policy in the following more interpretable form
\[
\phi_t^s(x) = K_t^s(x-(x^\star)_t^s),
\]
where $(x^\star)_t^s = -(K_t^s)^\dagger k_t^s$.
This has a convenient interpretation; to select $u_t$, we calculate the difference between $x$ and our \emph{desired} state $(x^\star)_t^s$ and then multiply that difference by the input gain matrix.
We can then interpret the policy as \emph{regulating} the state 
towards the desired state.


\section{Extended quadratic functions}
\label{sec:preliminaries}

In this section we describe extended quadratic functions, which are
quadratic functions with embedded linear equality constraints.
We explain how to verify attributes like convexity,
how they can be combined or pre-composed with an affine function,
and how to carry out partial minimization, where
we minimize over a subset of the variables.

An \emph{extended quadratic} function
$f:\reals^n\rightarrow\reals\cup\{+\infty\}$
has the form
\[
f(x)=
\frac{1}{2}\begin{bmatrix}x\\1\end{bmatrix}^T
\begin{bmatrix}P & q \\ q^T & r\end{bmatrix}
\begin{bmatrix}x \\ 1\end{bmatrix}
+ \mathcal I_{F,g} (x),
\]
where $P\in\symm^n$, $q\in\reals^n$, $r\in\reals$,
and $\mathcal I_{F,g}$ is the indicator function of 
the linear equality constraint $Fx+g =0$,
\[
\mathcal I_{F,g}(x) = \begin{cases}
    0 & Fx+g=0\\
    +\infty & \text{otherwise},
\end{cases}
\]
where $F\in\reals^{p \times n}$ and $g\in\reals^{p}$.
An extended quadratic function is specified by 
$P$, $q$, $r$, $F$, and $g$.
We refer to the function
\[
\frac{1}{2}\begin{bmatrix}x\\1\end{bmatrix}^T
\begin{bmatrix}P & q \\ q^T & r\end{bmatrix}
\begin{bmatrix}x \\ 1\end{bmatrix}
\]
as the \emph{quadratic part} of $f$, and
we refer to $\mathcal I_{F,g}(x)$ as the \emph{embedded equality constraints} in $f$.
We refer to $n$ as the \emph{dimension} of (the argument of) $f$.
We allow $p=0$, \ie, the case when there are no
embedded equality constraints in $f$.
In this case we refer to $f$ as a (non-extended) quadratic function.

\paragraph{Special cases and attributes of extended quadratic functions.}
An extended quadratic function $f$ is \emph{proper} 
\cite{rockafellar1970convex}
if there exists $x$ with $f(x) < +\infty$, \ie,
the embedded equality constraints are feasible.
An extended quadratic function $f$
is an \emph{extended quadratic form} 
if $q=0$, $r=0$, and $g=0$ (or $p=0$); in this case it is homogeneous 
of degree two.
If in addition there are no equality constraints,
\ie, $q=0$, $r=0$, and $p=0$, $f$ is a \emph{quadratic form}.
It is \emph{extended affine} if $P=0$, and \emph{affine} if in addition
there are no constraints.
An extended quadratic function $f$ is \emph{extended constant} if $P=0$ and
$q=0$.

\paragraph{Free parameter representation of equality constraints.}
The representation of the embedded 
equality constraints by $F$ and $g$ is evidently not unique.
For example, if $T \in \reals^{p \times p}$ is invertible,
$\tilde F = TF$ and $\tilde g=Tg$ gives another representation of the 
same constraints, \ie, $\mathcal I_{F,g} = \mathcal I_{\tilde F,\tilde g}$.
To resolve this non-uniqueness, and for other tasks as well,
it will be convenient to express the equality constraints in
\emph{free parameter form}, parametrized by $x_0 \in  \reals^n$
and $V_2 \in \reals^{n \times l}$, with $l = n-\Rank(F)$,
\BEQ\label{e-free-param}
\{x \mid Fx+g = 0\} = \{ V_2z+x_0 \mid z \in \reals^l \}.
\EEQ
Here $x_0$ is any particular solution of $Fx+g=0$, and
$\range(V_2) = \nullspace(F)$ (\ie, the nullspace of $F$).
Without loss of generality we can assume that $V_2^TV_2=I$.
(We will explain the subscript in $V_2$ shortly.)

We can determine whether the constraints are feasible, and if so,
find such a free parameter representation using the (full) 
singular value decomposition (SVD) of $F$, 
\[
F=\begin{bmatrix}U_1&U_2\end{bmatrix}
\begin{bmatrix}\Sigma & 0\\
0 & 0 \end{bmatrix}\begin{bmatrix}V_1 & V_2\end{bmatrix}^T,
\]
where $\Sigma \in \reals^{s \times s}$
contains the positive singular values of $F$, with $s = \Rank(F)$.
Then $Fx+g$ is feasible if and only if $U_2^T g =0$, and we can take
$x_0 = -V_1 \Sigma^{-1} U_1^T g = -F^\dagger g$, where $F^\dagger$ is the 
(Moore-Penrose) pseudo-inverse of $F$.
We can take $V_2$ as the matrix in our free parameter 
representation~(\ref{e-free-param}).
Finally, we note that we can replace the representation $F$ and $g$ with
$\tilde F = V_1^T$ and $\tilde g = \Sigma^{-1} U_1^T g$.
In this case $\tilde F$ satisfies $\tilde F \tilde F^T =I$, \ie,
its rows are orthonormal.
We refer to a representation of equality constraints with $FF^T=I$ (\ie, with 
orthonormal rows) as in \emph{reduced form}.
\emph{Reducing} an extended quadratic corresponds to converting its
equality constraints to reduced form, which can be done via the SVD as 
described above.

\paragraph{Equality of constraints.}
Using the decomposition above we can check whether two descriptions
of equality constraints (possibly of different dimensions $p$ and $\tilde p$)
are equal, \ie,  $\mathcal I_{F,g} = \mathcal I_{\tilde F,\tilde g}$.
Let $x_0,V_1,V_2$ and $\tilde x_0, \tilde V_1, \tilde V_2$ correspond
to the free parameter representation above
for $f$ and $\tilde f$ respectively.
Clearly we must have $\Rank(F)=\Rank(\tilde F)$, and in addition,
\[
V_1^T\tilde V_2 = 0, \quad V_1^T \tilde x_0 + g = 0, \quad
\tilde V_1^T V_2 = 0, \quad \tilde V_1^T x_0 + \tilde g = 0.
\]

\paragraph{Equality of extended quadratics.}
Two extended quadratics are equal, \ie, $f(x)=\tilde{f}(x)$ 
for all $x\in\reals^n$, if and only if $\mathcal I_{F,g} =
\mathcal I_{\tilde F,\tilde g}$ (discussed above),
and in addition
\[
f(x_0+V_2z)=\tilde{f}(x_0+V_2z) \quad \forall \; z \in\reals^r.
\]
Because $\mathcal I_{F,g} =
\mathcal I_{\tilde F,\tilde g}$, we can use the same free parameter representation.
This can be expressed as
\[
\begin{bmatrix}
V_2 & x_0 \\
0 & 1
\end{bmatrix}^T
\begin{bmatrix}
P & q \\
q^T & r
\end{bmatrix}
\begin{bmatrix}
V_2 & x_0 \\
0 & 1
\end{bmatrix}
=
\begin{bmatrix}
V_2 & x_0 \\
0 & 1
\end{bmatrix}^T
\begin{bmatrix}
\tilde P &\tilde  q \\
\tilde q^T &\tilde  r
\end{bmatrix}
\begin{bmatrix}
V_2 & x_0 \\
0 & 1
\end{bmatrix},
\]
or, more explicitly,
\BEAS
V_2^T P V_2 &=& V_2^T \tilde P V_2, \\
V_2^T Px_0 + V_2^T q &=&
V_2^T \tilde Px_0 + V_2^T \tilde q, \\
x_0^T P x_0 + 2x_0^T q + r &=&
x_0^T \tilde P x_0 + 2x_0^T \tilde q + \tilde r.
\EEAS

\paragraph{Convexity.}
Let $f$ be an extended quadratic function with free parameter representation $x_0$, $V_2$, as described above. We have that $f$ is convex if and only if
$V_2^T PV_2 \succeq 0$ ($A \succeq 0$ means the symmetric matrix 
$A$ is positive semidefinite).
It is strictly convex if and only if $V_2^T PV_2 \succ 0$ ($A \succ 0$ means $A$ is positive definite).

\begin{proof}
We prove the condition for convexity; the proof for strict convexity is almost identical.
Suppose $V_2^TPV_2 \succeq 0$.
Then the function $g(z)=f(x_0+V_2z)$ is convex in $z$, because
\[
\nabla^2_z g(z) = \nabla^2_z f(x_0+V_2z) = V_2^TPV_2 \succeq 0,
\]
where $\nabla_z^2 g(z)$ is the Hessian matrix of $g$.
It follows that $f$ is convex in $x$, since it is convex in the subspace of solutions to $Fx+g=0$ and $+\infty$ otherwise.

Suppose that there exists $v\in\nullspace(F)$ such that $v^TPv < 0$, and that $f$ is convex.
Then $f$ is not convex along the line $x_0+tv$ for $t\in\reals$, because
\[
\frac{\partial^2}{\partial t^2}f(x_0+tv) = v^TPv < 0.
\]
It follows that $f$ is not convex, which is a contradiction.
\end{proof}

\paragraph{Nonnegativity.}
An extended quadratic is nonnegative if and only if
\[
\begin{bmatrix}
V_2 & x_0 \\
0 & 1
\end{bmatrix}^T
\begin{bmatrix}
P & q \\
q^T & r
\end{bmatrix}
\begin{bmatrix}
V_2 & x_0 \\
0 & 1
\end{bmatrix}
\succeq 0.
\]

\paragraph{Sum.}
The sum of two extended quadratics of the same dimension,
\[
f(x) + g(x) = \frac{1}{2}\begin{bmatrix}x\\1\end{bmatrix}^T
\begin{bmatrix}P & q \\ q^T & 
r\end{bmatrix}
\begin{bmatrix}x \\ 1\end{bmatrix}
+
\mathcal I_{F,g}(x)
+
\frac{1}{2}\begin{bmatrix}x\\1\end{bmatrix}^T
\begin{bmatrix}\tilde{P} & \tilde{q} \\ \tilde{q}^T & 
\tilde{r}\end{bmatrix}
\begin{bmatrix}x \\ 1\end{bmatrix}
+
\mathcal I_{\tilde F, \tilde g}(x)
\]
is also an extended quadratic, 
with quadratic part
\[
\frac{1}{2}\begin{bmatrix}x\\1\end{bmatrix}^T
\begin{bmatrix}P+\tilde{P} & q+\tilde{q} \\ (q+\tilde{q})^T & 
r+\tilde{r}\end{bmatrix}
\begin{bmatrix}x \\ 1\end{bmatrix},
\]
and equality constraints
\[
    \begin{bmatrix}F \\ \tilde F\end{bmatrix}
    x  + \begin{bmatrix} g \\ \tilde g \end{bmatrix} = 0.
\]
The sum can be improper, even when $f$ and $g$ are not.
After adding two extended quadratics, we can reduce the equality constraints
(which also checks if the sum is proper).

\paragraph{Scalar multiplication.}
We \emph{define}
scalar multiplication of an extended quadratic, $\alpha f$, as
the extended quadratic with the same equality constraints,
and quadratic part scaled by $\alpha$:
\[
(\alpha f)(x)
= 
\frac{1}{2}\begin{bmatrix}x\\1\end{bmatrix}^T
\begin{bmatrix}\alpha P & \alpha q \\ \alpha q^T & \alpha r\end{bmatrix}
\begin{bmatrix}x \\ 1\end{bmatrix}
+ \mathcal I_{F,g}(x).
\]
If $f$ is convex and $\alpha \geq 0$, $\alpha f$ is convex.
Not that when $\alpha<0$, our definition of $\alpha f$ is not the usual
mathematical one, since ours takes the value $+\infty$ when the implicit
constraints are violated, whereas the under the usual definition,
$\alpha f$ would take the value $-\infty$.

\paragraph{Affine pre-composition.}
Suppose that $x = h(z)= Az+b$ is an affine function of $z$, and 
consider $g = f\circ h$, \ie, $g(z) = f(Az+b)$.
The function $g$ is extended quadratic, and has the form
\[
g(x) =
\begin{bmatrix}z \\ 1\end{bmatrix}^T
\begin{bmatrix}\tilde P & \tilde q \\ \tilde q^T & \tilde r\end{bmatrix}
\begin{bmatrix}z \\ 1\end{bmatrix} + 
\mathcal I_{\tilde F,\tilde g}(z),
\]
where 
\[
\begin{bmatrix}\tilde P & \tilde q \\ \tilde q^T & \tilde r\end{bmatrix} =
\begin{bmatrix}
        A & b \\
        0 & 1
    \end{bmatrix}^T
\begin{bmatrix} P & q \\ q^T & r\end{bmatrix}
\begin{bmatrix}
        A & b \\
        0 & 1
    \end{bmatrix},
\qquad
\tilde F = FA, \qquad \tilde g = Fb+g.
\]
The equality constraints can be reduced.
If $f$ is convex, $g$ is convex.

\paragraph{Partial minimization.}

Next we consider partial minimization of an extended quadratic function, meaning we fix
a subset of its variables and minimize over the other variables.
There are two cases to consider:
\begin{itemize}
	\item \emph{Strictly convex.} When the function is strictly convex in the variables we are minimizing over, the function of the remaining variables is always extended quadratic.
	\item \emph{Convex but not strictly convex.} When the function is convex (but not strictly convex) in the variables we are minimizing over, and a technical condition holds (which always holds in the strictly convex case), the function of the remaining variables is always extended quadratic.
\end{itemize}

Suppose $f$ is an extended quadratic function of two variables $x\in\reals^n$ and $u\in\reals^m$, \ie, it has the form
\[
f(x,u) = 
\frac{1}{2}\begin{bmatrix}x\\u\\1\end{bmatrix}^T
\begin{bmatrix}
P_{xx} & P_{xu} & q_x \\
P_{ux} & P_{uu} & q_u \\
q_x^T & q_u^T &  r
\end{bmatrix}
\begin{bmatrix}x \\u\\ 1\end{bmatrix}
+
\mathcal{I}_{F,g}(x,u),
\]
where $F=\begin{bmatrix}F_x & F_u\end{bmatrix}$.
Also, suppose that $f$ is convex in $u$ for all $x$.
The function
\[
g(x) = \inf_u f(x,u)
\]
gives the partial minimization of $f$ over $u$.

Evaluating $\inf_u f(x,u)$ is itself a convex optimization problem, and for it to be feasible, the (extended quadratic) function
$h(u)=f(x,u)$
must be proper.
We have that $h$ is proper if and only if
\[
x \in \{x \mid F_x x + g \in \range(F_u)\}
= \{x \mid \tilde{F}x + \tilde{g} = 0\},
\]
where
$\tilde F = (I-F_uF_u^\dagger)F_x$ and $\tilde g = (I-F_uF_u^\dagger)$,
which is a linear equality constraint on $x$.
We can express $g$ in the equivalent form
\BEQ
g(x)=\inf_uf(x,u) +
\mathcal{I}_{\tilde{F},\tilde{g}}(x).
\label{eq:g}
\EEQ

We can convert the equality constraint on $x$ above into its free parameter representation, parameterized by $x_0\in\reals^n$ and $V_2\in\reals^{n\times l}$, \ie,
\[
\{x \mid \tilde{F}x + \tilde{g} = 0\}=\{V_2z + x_0 \mid z \in\reals^l\}.
\]


Carrying out the partial minimization to find $g$
amounts to solving the (equality constrained quadratic) optimization problem
\[
\begin{aligned}
& \text{minimize}
& &
\frac{1}{2}\begin{bmatrix}u\\1\end{bmatrix}^T
\begin{bmatrix}
P' & q'  \\
q'^T & r'  
\end{bmatrix}
\begin{bmatrix}u\\ 1\end{bmatrix}
\\
& \text{subject to}
&& F_{u} u + g'=0,
\end{aligned}
\]
where
\[
P'=P_{uu}, \quad q' = P_{ux}x + q_u, \quad g' = F_x x + g, \quad r'=x^T P_{xx} x + x^T q_x + q_x^T x + r.
\]
We replace $x$ with its free parameter representation $x_0+V_2z$,
so that the optimization problem is guaranteed to be feasible.
The KKT conditions for this problem for $x^\star$ and $\nu^\star$ to be optimal~\cite{boyd2004convex} are
\[
\begin{bmatrix}
P_{uu} & F_u^T \\
F_u & 0
\end{bmatrix}
\begin{bmatrix}u^\star\\\nu^\star\end{bmatrix}
=
-
\begin{bmatrix}
q'\\g'
\end{bmatrix}
=
-\begin{bmatrix}P_{ux}\\ F_x\end{bmatrix}V_2z - \begin{bmatrix}P_{ux}x_0 + q_u \\ F_x x_0 + g\end{bmatrix}.
\]
This linear system has a solution for all $z$ if and only if
\BEQ\range(\begin{bmatrix}P_{uu}&F_u^T\\F_u&0\end{bmatrix})\supseteq\range(\begin{bmatrix}P_{ux}V_2 & P_{ux}x_0+q_u\\F_xV_2 & F_xx_0+g\end{bmatrix})
.
\label{eq:range}
\EEQ
This is the technical condition that we have been referring to.
This condition is guaranteed to hold if, \eg, $f$ is strictly convex in $u$.
If the technical condition \eqref{eq:range} holds,
then we can express a $u^\star$ as
\[
u^\star =-\begin{bmatrix}I & 0\end{bmatrix}\begin{bmatrix}
P_{uu} & F_u^T \\
F_u & 0
\end{bmatrix}^\dagger
\left(
\begin{bmatrix}P_{ux}\\ F_x\end{bmatrix}x + \begin{bmatrix}q_u\\g\end{bmatrix}
\right)
= Ax+b,
\]
where $A\in\reals^{m\times n}$ and $b\in\reals^m$.
This always satisfies the constraint $F_xx+F_uu^\star+g=0$.

Plugging this $u^\star$ back into \eqref{eq:g}, we find that $g$ is an extended quadratic, with quadratic part
\BEAS
f(x,u^\star)&=&
\frac{1}{2}
\begin{bmatrix}x\\1\end{bmatrix}^T
\begin{bmatrix}I&0\\A&b\\0&1\end{bmatrix}^T
\begin{bmatrix}
P_{xx}&P_{xu}&q_x\\
P_{ux}&P_{uu}&q_u\\
q_x^T&q_u^T&1
\end{bmatrix}
\begin{bmatrix}I&0\\A&b\\0&1\end{bmatrix}
\begin{bmatrix}x\\1\end{bmatrix},
\\
&=&\frac{1}{2}
\begin{bmatrix}x\\1\end{bmatrix}^T
\begin{bmatrix}
P_{xx}+P_{xu}A+A^TP_{ux}+A^TP_{uu}A&(P_{xu}b+A^TP_{uu}b+A^Tq_u+q_x)\\
(P_{xu}b+A^TP_{uu}b+A^Tq_u+q_x)^T&b^TP_{uu}b + b^Tq_u+q_u^Tb+r
\end{bmatrix}
\begin{bmatrix}x\\
1\end{bmatrix},
\EEAS
and equality constraints
$\mathcal{I}_{\tilde{F},\tilde{g}}(x)$.

If \eqref{eq:range} does not hold, or $f$ is nonconvex in $u$, then there is at least one $x$ where $g(x)=-\infty$ and $g$ is no longer an extended quadratic function of $x$.
We can check the technical condition \eqref{eq:range} by noting that
\[
\range(A) \supseteq \range(B) \iff (I-AA^\dagger)B = 0.
\]

\section{Dynamic programming solution}

In this section we use dynamic programming to show that the solutions
to the problems we consider in this paper have the form \eqref{eq:policies}.

\subsection{Finite-horizon}

The cost-to-go (or value) function $V_t^s: \reals^n \rightarrow \reals \cup \{-\infty,+\infty\}$ is defined as the cost
achieved by an optimal policy starting at time $t$ from a given state and mode, or
\[
V_t^s(x) =
\inf_{\phi_t}
\Expect \left[\sum_{\tau=t}^{T-1} g_\tau^{s_\tau}(x_\tau,\phi_\tau^{s_\tau}(x_\tau),w_\tau) + g_T^{s_T}(x_T,w_T) \right], \quad t=0,1,\ldots,T,
\]
subject to the dynamics \eqref{eq:dynamics}, $x_t=x$, and $s_t=s$.
Given a collection of functions $V=(V^1,\ldots,V^K)$, define the Bellman operator $\mathcal{T}_t$
applied to that collection as
\[
(\mathcal{T}_t V)^s(x) = \inf_u
\Expect_{w_t} \Expect_{s'}
\left[
g_t^s(x,u,w_t)+ V^{s'}(f_t^s(x,u,w_t))
\right]
\]
where $s'=i$ with probability $\Pi_{t,ij}$ if $s=j$.
It is well known that the cost-to-go functions $V_0,V_1,\ldots,V_T$ satisfy the dynamic programming
(DP) recursion~\cite{bertsekas2017dynamic,bertsekas2012dynamic},
\[
V_t = \mathcal{T}_t V_{t+1}, \quad t=T-1,\ldots,0,
\]
with
$V_T^s(x) = g_T^s(x)$.

Defining the state-action cost-to-go function
$Q_t^s: \reals^n \times \reals^m \rightarrow \reals \cup \{-\infty,+\infty\}$ as
\[
Q_t^s(x,u)=\Expect_{w_t} \Expect_{s'}
\left[
g_t^s(x,u,w_t)+ V_{t+1}^{s'}(f_t^s(x,u,w_t))
\right],
\]
the optimal policies are given by
\[
\phi_t^s(x) = \argmin_u Q_t^s(x,u).
\]

\paragraph{Extended quadratic cost-to-go functions.}
We show that, barring pathologies, the cost-to-go functions
$V_t^s$ are extended quadratic functions of $x$ for $t=0,\ldots,T$.
Intuitively, this is because the Bellman operator preserves
the `extended quadraticity' of the cost-to-go functions.

We show this by induction.
The last cost-to-go function $V_T^s$ is extended quadratic in $x$ by definition.
Suppose, then, that $V_{t+1}^s$ is extended quadratic in $x$.
Then $Q_t^s$ must be extended quadratic in $x$ and $u$ because
it is the expectation of an extended quadratic
plus an extended quadratic composed with an affine function.
Because $V_t$ is equal to the partial minimization of
$Q_t^s$, an extended quadratic, $V_t^s$ will be an extended quadratic function of $x$.
(If $V_t^s$ is not an extended quadratic function of $x$, then the cost is either $+\infty$ or $-\infty$, a pathology.)
Therefore, the
cost-to-go functions are extended quadratic functions of $x$.

\paragraph{Affine policies.}
Because $\phi_t^s$ is equal to the solution of partially minimizing an extended quadratic function $Q_t^s$, it follows that
there exists an optimal policy that is affine in $x$ for each $t,s$, and has the form in~\eqref{eq:policies}.

\subsection{Infinite-horizon}
\label{sec:infinite}

The cost-to-go function of the infinite-horizon problem,
$V: \reals^n \rightarrow \reals \cup \{-\infty,+\infty\}$, is given by
\[
V^s(x) = \inf_\phi \Expect \left[
\sum_{t=0}^\infty \gamma^t g^{s_t}(x_t,\phi^{s_t}(x_t),w_t)
\right],
\]
subject to the dynamics~\eqref{eq:dynamics}, $x_0=x$, $s_0=s$, where
the expectation is over $w_t$.
Given a collection of functions $H=(H^1,\ldots,H^K)$ for $H^s:\reals^n \rightarrow \reals$,
define the Bellman operator $\mathcal{T}$ applied to the collection as
\[
(\mathcal{T} H)^s(x) = \inf_u
\Expect_{w} \Expect_{s'}
\left[
g^s(x,u,w) + \gamma H^{s'}(f^s(x,u,w))
\right]
\]
where $s'=i$ with probability $\Pi_{ij}$ if $s=j$.
It is well known that the cost-to-go function is the unique fixed point
of the Bellman operator~\cite{bertsekas2017dynamic,bertsekas2012dynamic}, or
\[
V = \mathcal{T} V,
\]
and that
\[
V = \lim_{k\rightarrow \infty} \mathcal{T}^k V_0,
\]
for any bounded function $V_0^s: \reals^n \rightarrow \reals$.

Defining the state-action cost-to-go function $Q^s: \reals^n \times \reals^m \rightarrow \reals \cup \{-\infty,+\infty\}$ as
\[
Q^s(x,u) = \Expect_{w} \Expect_{s'}
\left[
g^s(x,u,w) + \gamma V^{s'}(f^s(x,u,w))
\right]
\]
where $s'=i$ with probability $\Pi_{ij}$ if $s=j$,
the optimal policy is given by
\[
\phi^s(x) = \argmin_u Q^s(x,u).
\]

\paragraph{Quadratic cost-to-go function.}

We show that the cost-to-go function is an extended quadratic function of $x$.
If $H$ is extended quadratic in $x$, so is $\mathcal{T} H$ by the same logic as in the finite-horizon case (barring pathologies).
Starting with $V_0(x)=0$, a bounded extended quadratic function,
we have that $\mathcal{T}^k V_0$ is an extended quadratic function of $x$ for $k\in\mathbb{N}$.
Therefore, its limit, the cost-to-go function $V$, is an extended quadratic function of $x$.

\paragraph{Affine policy.}
We have that $Q$ is extended quadratic because $V$ is extended quadratic.
Therefore, by the same logic as the finite-horizon case, there exists an optimal policy that is affine in $x$.

\subsection{Avoiding pathologies}

Pathologies are most likely to manifest when one performs partial minimization of $Q_t$ to find $V_t$.
If $Q_t$ is strictly convex in the variables one is minimizing over, the partial minimization will always yield an extended quadratic function of the other variables.
We can enforce this by making $g_t^s(x,u)$ convex in $x$ and strictly convex in $u$, but pathologies can still occur, \eg, the cost could be infinite.
(This is done in LQR, where $g(x,u)=x^TQx+u^TRu$ and $R \succ 0$.)

If $Q_t$ is convex but not strictly convex, then in addition, the range condition \eqref{eq:range} must hold.
The best way to check this is to solve the problem and see (numerically) if the range condition holds.
If it does not, then the problem is ill-posed and the cost or dynamics should be changed.

\subsection{Linear policies}

In some problems, \eg, in classical LQR, the policies are linear functions of $x$.
In this section we come up with a sufficient condition on problems that have linear (optimal) policies.

If the dynamics are linear,
and the stage cost is a convex quadratic form with homogeneous equality
constraints, then
the optimal policies will be linear.
This is because the cost-to-go functions are quadratic forms with homogeneous equality constraints, \ie, they are of the form
\[
V_t^s(x) = \frac{1}{2}\begin{bmatrix}x\\1\end{bmatrix}\begin{bmatrix}P&0\\0&r\end{bmatrix} + \begin{cases}0 & Fx=0 \\ +\infty & \; \text{otherwise}\end{cases}.
\]
This follows because quadratic forms with homogeneous equality constraints are closed under the same operations as extended quadratics.
It is easy to show that the solution of partially minimizing these functions is a linear function of the other variables, meaning the policies will be linear.

\section{Implementation}

In this section we describe how to numerically perform dynamic programming to
find the cost-to-go functions, state-action cost-to-go functions, and policies, in the finite and infinite-horizon problems.
One of the main difficulties in carrying out dynamic programming in the general case is representing the functions of interest.
For the problems that we consider, representation is easy, since
we can represent the cost-to-go functions, state-action cost-to-go functions, 
and policies via explicit lookup tables, indexed by time and mode, of the coefficients
in the associated extended quadratic (or affine) functions.

To perform dynamic programming, all we need to do is apply the Bellman operator.
The Bellman operator requires several operations on extended quadratic functions: addition, scalar multiplication (in the infinite-horizon case), affine composition, partial minimization, and expectation.
The first four can be carried out using elementary linear algebra operations, as described in \S\ref{sec:preliminaries}.
The remaining operation is expectation.

\paragraph{Expectation.}
There are two expectations: one over the next mode and one over the disturbance.
The expectation over the next mode is easy, since there are $K$ possibilities which we can enumerate by replacing the expectation by a weighted sum.
This leaves the expectation over $w_t$.

If we only have access to a sampling oracle, we can approximate this expectation using Monte Carlo expectation, \ie, we sample $w_1,\ldots,w_N$ for some large number $N$, calculate the corresponding extended quadratic function for each $w_i$, and then average those extended quadratic functions.
This is, in a way, the best that we can do given only an oracle that provides independent samples of the problem data.

When the cost-to-go function is a (non-extended) quadratic function, \ie, it has no equality constraints, we can exactly perform expectation if we know the first and second moments of the dynamics matrices and the first moment of the cost matrix, and nothing more.
(See Appendix~\ref{sec:exact} for more on this.)
If, however, the cost-to-go function contains equality constraints, then we need more knowledge than just the first and second moments of the dynamics matrices, and we need to fall
back to approximating the expectation using the Monte Carlo procedure explained above.

In our implementation we use Monte Carlo expectation because it requires nothing more than 
a sampling oracle, and is easy to implement.
The full version of the finite-horizon algorithm with Monte Carlo expectations is presented in Algorithm~\ref{alg1}.
\begin{algdesc}
\emph{Finite-horizon extended quadratic control.} 
\begin{tabbing}
{\bf given} $T$, $N$, $K$, $A_t^s(w)$, $B_t^s(w)$, $c_t^s(w)$, $g_t^s(w)$, $\Pi_t^s$, $g_T^s$, and independent samples $w_1^{t,s},\ldots,w_N^{t,s}$.\\
Set $V_T^s = g_T^s$.\\
{\bf for} $t=T-1,\ldots,0$\\
\qquad \= 1.\ $Q_t^s(x,u) = \frac{1}{N}\sum_{i=1}^N\sum_{s'=1}^K \Pi_{s',s}^t \left(g_t^s(x,u,w_i^{t,s}) + V_{t+1}^{s'}(A_t^s(w_i^{t,s}) x + B_t^s(w_i^{t,s}) u + c_t^s(w_i^{t,s}))\right)$.\\
\> 2.\ \emph{Partial minimization.} Form $V_t^s(x)=\inf_u Q_t^s(x,u)$ and $\phi_t^s(x) = K_t^s x + k_t^s$.\\
{\bf end for}.
\end{tabbing}
\label{alg1}
\end{algdesc}

\paragraph{Infinite-horizon.}

To perform infinite-horizon dynamic programming, we simply call finite-horizon dynamic programming with a final stage cost function of zero, with stage cost functions multiplied by $\gamma^t$, and with the time horizon being the number of times to apply the Bellman operator.
One could devise a more sophisticated termination condition, \eg, terminating when $V_t$ is ``close'' to $V_{t+1}$, rather than applying the Bellman operator a fixed number of times.
However, in practice, the cost-to-go function converges after a small number (\eg, $10$ to $20$) of iterations.

\paragraph{Python implementation.}

We have developed an open-source Python library that implements Algorithm~\ref{alg1}.
The main object is \verb+ExtendedQuadratic+, which is initialized by supplying $P$, $q$, $r$, $F$, and $g$.
One can perform arithmetic operations on \verb+ExtendedQuadratic+s,
as well as perform
affine pre-composition, partial minimization, and check
equality, convexity, and nonnegativity.
We also provide a function that maps \verb+ExtendedQuadratic+s to cvxpy~\cite{diamond2016cvxpy} expressions and vice versa.

The main methods are \verb+dp_finite+ and \verb+dp_infinite+, which implement
finite and infinite horizon dynamic programming, respectively.
One supplies a sampling function \verb+sample(t,N)+, which provides a batch sample problem data (of size $N$) for time $t$.
(For \verb+dp_infinite+, the sample function does not take $t$ as an argument as it should be time invariant.)
The functions take two additional arguments: the number of Monte Carlo samples, the time horizon (in the infinite horizon function, this is the number of times to apply the Bellman operator), and the discount factor (in the infinite-horizon case).
The method returns the cost-to-go functions $V_t^s(x)$, the state-action cost-to-go functions $Q_t^s(x,u)$, and the (optimal) policies $(K_t^s,k_t^s)$.

\paragraph{Runtime.}
The finite-horizon extended quadratic control algorithm, \ie, $T$ applications of the Bellman operator, requires approximately
\[
TK^2N\max\{1,p\}\max\{n,m\}^2
\]
operations.
Our (na\"ive) single-threaded Python implementation applied to a random
moderately sized problem with $n=25$, $m=50$, $N=100$
(number of Monte Carlo samples), $K=5$,
and $T=25$ Bellman operator evaluations
takes about $12.8$ seconds to calculate the optimal policies
on a six-core $3.7$ GHz Intel i7.
The bulk of the computational effort lies in calculating the Monte Carlo expectation, which one could parallelize across multiple CPUs or GPUs to make the algorithm faster (see the MPI implementation below).

\paragraph{MPI Implementation.}

We mentioned above that the algorithm could be significantly sped up with a parallel implementation.
We have developed a distributed implementation using the Message passing interface (MPI)~\cite{dongarra1993proposal}.
MPI is a language-independent message-passing standard designed for parallel computing, and is the dominant model in high-performance computing applications today.
Though there are multiple ways to implement these algorithms in MPI, perhaps the simplest way is to parallelize the (Monte Carlo) sum over $i$ in line $1$ of Algorithm~\ref{alg1}.
If we have $r$ processors, we can set $N$ to a multiple of $r$ and have each processor perform a (smaller) Monte Carlo expectation over $N/r$ samples, and then \emph{reduce} these by averaging them.
This can provide significant reductions in runtime, provided $N$ is substantially greater than the number of processors.

\paragraph{Measuring Monte Carlo error.}

Our calculation of the cost-to-go functions and policies is approximate, because we use a Monte Carlo expectation instead of an exact expectation.
We can measure the error in our procedure by running it multiple times with different random seeds and checking how much the cost-to-go functions and policies vary across the runs.
We could also use this idea to \emph{dynamically} select the number of samples used in the Monte Carlo expectations to get a solution that is within a prescribed error.

\section{Extensions and variations}
\label{sec:extandvars}

There are several extensions of and variations on
the problems that we describe in this paper.

\subsection{Tractable extensions and variations}

In this section, we describe extensions and variations that maintain tractability, meaning we can still efficiently solve for the optimal policy with
a slightly modified version of the algorithms described in this paper.

\paragraph{Arbitrary information patterns.}

We can extend the algorithms described in this paper to problems with arbitrary information patterns.
One example of a different information pattern is where we have access to part of the disturbance before selecting $u_t$.
Suppose that the disturbance is partitioned into $w^b$ (before) and $w^a$ (after), \ie, $w=(w^b,w^a)$, and that we get to see $w^b$ before we select the input.
Our policy then becomes a function of $w^b$ and we arrive at the following modified dynamic programming recursion (we work with the infinite-horizon case for ease of notation):
\BEAS
\phi^s(x,w^b) &=& \inf_u
\Expect_{w^a \mid w^b}
\Expect_{s'}
\left[
g^s(x,u,(w^b,w^a)) + \gamma V^{s'}(f^s(x,u,(w^b,w^a)))
\right]
\\
V^s(x) &=& \inf_{\phi}
\Expect_{w=(w^b,w^a)}
\Expect_{s'}
\left[
g^s(x,\phi^s(x,w^b),w) + \gamma V^{s'}(f^s(x,\phi^s(x,w^b),w))
\right].
\EEAS
In the absence of pathologies, the cost-to-go functions are still extended quadratic functions of $x$, and the policies are still affine functions of $x$, though the notation is substantially more complicated and the algorithm is slower because of the nested Monte Carlo iterations.
Also, evaluating the policy requires a Monte Carlo expectation followed by solving an equality constrained quadratic optimization problem.

\paragraph{Incorporating a finite action.}

We can incorporate a finite action $a_t\in\{1,\ldots,A\}$ that we select at each time step that alters
the state dynamics, mode switching dynamics, and cost.
If we can select the finite action based on the current mode \emph{and} state, the problem quickly becomes intractable.
However, if we enforce that the action can only be selected based on the mode (and not $x_t$), meaning
\[
a_t = \psi(s_t),
\]
where $\psi:\{1,\ldots,K\} \rightarrow \{1,\ldots,A\}$ is called the \emph{action policy}.

Under this information pattern, we show that, for every initial state, there is an optimal policy that is affine.
We provide a sketch of the proof in the infinite-horizon problem.
Suppose the initial state $x_0=x^\text{init}$.
There are a finite number (in fact, $K^A$) of action policies (unfortunately, there are many more in finite-horizon problems).
Suppose we fix an action policy.
This results in an extended quadratic control problem, because the state dynamics, mode switching dynamics, and cost depend only on $s$.
This means that there exists an optimal state policy that is affine.
So we can construct an optimal policy by, for every action policy, calculating the average cost starting from $x_0$ with an optimal state policy, and then taking the action policy (and its corresponding affine control policy) that achieves the lowest cost as optimal.
Thus we have shown that there is an optimal policy that is affine.

To the best of our knowledge, there is no way to find the action policy faster than brute force (in the worst case).
But for a small problem, \eg, $K=5$ and $A=4$, we would only have to solve $625$ extended quadratic problems to find the \emph{optimal} action policy and state policy.

\paragraph{Mode switching dependent on the disturbance.}

We mentioned in \S\ref{sec:introduction} that we assumed $\Pi$ was independent of
$w$.
We assumed this for simplicity, but it is technically not required.
One could make the switching probability matrix depend on the disturbance, and this would not affect the Bellman
operator, except for the fact that one would also have to also provide samples of $\Pi$, because it becomes random.

\paragraph{Average cost problems.}

In the \emph{average cost problem}, we assume that the dynamics, cost,
and policy are time-invariant, and our goal is to find $\phi$ that
minimizes the expected average stage cost, or
\[
\lim_{T\rightarrow\infty}
\frac{1}{T+1}
\Expect
\left[
\sum_{t=0}^T g^{s_t}(x_t,\phi^{s_t}(x_t),w_t)
\right].
\]
Average cost problems with the same dynamics and cost as the problems considered in this paper may be solved using slightly modified versions of the algorithms described in this paper.

\paragraph{Periodic problems.}

Periodic problems are infinite-horizon problems with \emph{periodic} dynamics and cost (see the paragraph ``Periodic Problems'' in Chapter 4.6 of \cite{bertsekas2012dynamic}).
Periodic problems with the same dynamics and cost as the problems considered in this paper may be solved using slightly modified versions of the algorithms described in this paper.

\subsection{Heuristic extensions and variations}

In this section we describe extensions and variations that are heuristic, meaning we search for a suboptimal policy with no guarantee of optimality.
Though heuristic, they can be very effective in practice.

\paragraph{Nonlinear systems.}

We can extend the algorithms described in this paper to approximately solve stochastic control problems with nonlinear dynamics and extended quadratic stage costs.
(It is not hard to extend this to non-quadratic cost and infinite horizon.)
Suppose that we have time-invariant dynamics
\[
x_{t+1}=f(x_t,u_t,w_t)
\]
where $f$ is nonlinear and differentiable in $x_t$ and $u_t$.
Suppose that our system starts in the state $x_0=x^\text{init}$ and that at iteration $k$ of the algorithm we have fixed a policy $u_t=\phi_k(x_t)$.

At each time step $t$, we can approximate the nonlinear dynamics with the following (random, time-varying) affine dynamics
\[
x_{t+1} \approx A_t(w_t)x_t + B_t(w_t)u_t+c_t(w_t),
\]
where
\BEQ
\label{aab}
\begin{aligned}
A_t(w_t) &:= D_x f(x_t,u_t,w_t)\\
B_t(w_t) &:= D_u f(x_t,u_t,w_t)\\
c_t(w_t) &:= f(x_t,u_t,w_t) - D_xf(x_t,u_t,w_t)x_t - D_u f(x_t,u_t,w_t) u_t.
\end{aligned}
\EEQ
Here $D_x f(x,u)$ is the Jacobian matrix of $f$ at $x$.

To generate a sample of $A_t$, $B_t$, and $c_t$, we first generate a sample trajectory using the current policy by repeatedly applying the state transition function, \ie,
\[
x_{t+1} = f(x_t,\phi_k(x_t),w_t),
\]
starting with $x_0=x^\text{init}$,
and plugging the resulting states, controls, and disturbances into \eqref{aab}.

We can now solve the stochastic problem assuming that instead we have the time-varying affine dynamics described above, resulting in an affine control policy $\phi_{k+1}$ which we use in the next iteration of the algorithm.
There are no guarantees on convergence or optimality of this method, but 
it seems like a natural generalization 
of iLQR/iLQG~\cite{li2004iterative,todorov2005generalized} to general stochastic 
nonlinear systems, where one allows noise that is multiplicative in the state and control.



\paragraph{Convex constraints on the input.}

A natural question to ask is: how does one incorporate convex constraints
on the input?
For example, we might want to constrain the input to be within some range, \eg, $u^\text{min}\leq u_t \leq u^\text{max}$.
Incorporating such constraints directly into dynamic programming makes it intractable.
Even though the cost-to-go functions are convex, it is impossible in general to compute them, let alone represent them.

However, as an approximation, we can solve the problem that ignores the constraints on the
input, which we can do exactly, and then incorporate
the constraints when we actually select $u_t$.
We can then incorporate the convex constraints with the following ADP (or control-Lyapunov) policy
\BEQ
\phi_t^s(x) = \inf_{u \in U} Q_t^s(x,u),
\label{eq:adp}
\EEQ
which calculates the optimal policy using $Q_t^s$ as the \emph{approximate} state-action cost-to-go function.
Because \eqref{eq:adp} is a convex optimization problem, we can exactly (and efficiently) solve it.
Although this is an approximation, we are guaranteed to select inputs that satisfy the constraints.

\section{Applications}

In this section we describe several known and new applications.

\subsection{LQR}

LQR is a classical problem in control theory, first identified and solved by Kalman in the late 1950's~\cite{kalman1960contributions}.
There are many variations of LQR;
in this example we focus on the infinite horizon LQR problem.
The system has a time-invariant linear state transition function, meaning its dynamics are described by
\[
x_{t+1} = Ax_t + Bx_t + w_t, \quad t=0,1,\ldots,
\]
where $\Expect w_t=0$ and $\Expect[w_tw_t^T] = W$,
and our stage cost is a quadratic form, \ie,
\[
g(x,u) = x^TQx+u^TRu,
\]
where $Q\succeq 0$ and $R \succ 0$.

It is well known that the cost-to-go function is of the form
$V(z) = z^TPz$
where $P\succeq 0$ satisfies the algebraic Riccati equation (ARE)
\[
P = Q + A^TPA - A^TPB(R+B^TPB)^{-1}B^TPA,
\]
and that the optimal policy is linear state feedback $u_t=Kx_t$, where
\[
K = -(R+B^TPB)^{-1}B^TPA.
\]

Infinite-horizon LQR (and all of its tractable variants) are instances of the problem that we describe in this paper, because the dynamics are affine and the cost is a quadratic form.
However, there is only one mode, the cost, dynamics, and input matrices do not depend on $w_t$,
and the stage cost is a convex quadratic form, not an extended quadratic.
We know, without deriving the ARE,
that the cost-to-go function is a (non-extended) quadratic form
and that the optimal policy is linear, and we can efficiently calculate them
using the algorithms described in this paper.
It is worth noting that there are many other specialized (and efficient) methods of exactly solving the ARE; see, \eg,~\cite{laub1979schur,anderson2007optimal}.
LQR serves as a good test of our numerical implementation, since
we can compare the cost-to-go functions and policies that we find with specialized solvers 
for the ARE.

\subsection{LQR with random dynamics}

We can easily extend infinite-horizon LQR to incorporate random dynamics and input matrices.
That is, our dynamics are described by
\[
x_{t+1} = A_t x_t + B_t u_t + c_t, t=0,1,\ldots,
\]
where $A_t$, $B_t$, and $c_t$ are (jointly) random, and we have a quadratic form stage cost.
$A_t$, $B_t$, and $c_t$ can have \emph{any} joint distribution.

This problem was first identified and solved by Drenick and Shaw in 1964~\cite{drenick1964optimal} and then in continuous time by Wonham in 1970~\cite{wonham1970random}.
For a more modern treatment in discrete time, see the paragraph ``Random System Matrices'' on pages 123-124 in~\cite{bertsekas2017dynamic}.
The cost-to-go function in infinite-horizon LQR with random system matrices can diverge to $+\infty$ if there is too much noise in the system; this is referred to as the uncertainty threshold principle~\cite{athans1977uncertainty}, and is a great example of a pathology.

\paragraph{Numerical example.}

We reproduce the results in the original paper on the uncertainty threshold principle~\cite{athans1977uncertainty}.
Here we have a one-dimensional system ($n=m=1$) with dynamics described by
\[
x_{t+1} = a x_t + b u_t,
\]
where $a\sim\mathcal{N}(\bar a, \Sigma_{aa})$ and $b\sim\mathcal{N}(\bar b, \Sigma_{bb})$ (where $\mathcal{N}(\mu,\sigma)$ is the normal distribution with mean $\mu$ and standard deviation $\sigma$), and our stage cost is
\[
g(x,u) = x^2 + u^2.
\]
    


We solve the finite-horizon problem (with zero final cost), resulting in  cost-to-go functions that have the form
\[
V_t(x) = k_t x^2, \quad t=0,1,\ldots,T.
\]

\begin{figure}
    \centering
    \begin{subfigure}{.4\textwidth}
	    \includegraphics[width=\linewidth]{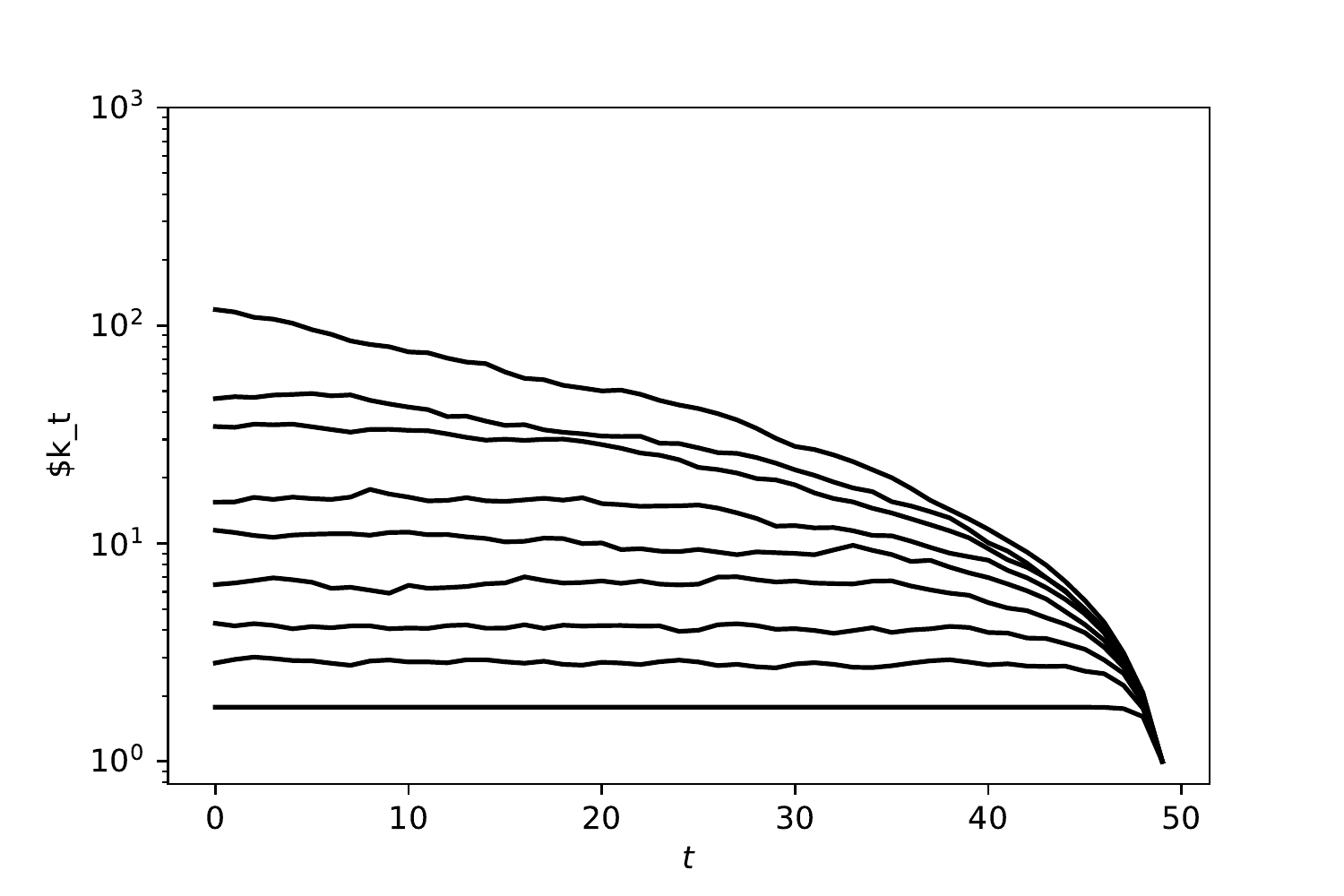}
	    \caption{}
	    \label{fig:random1}
    \end{subfigure}
    \begin{subfigure}{.4\textwidth}
	    \includegraphics[width=\linewidth]{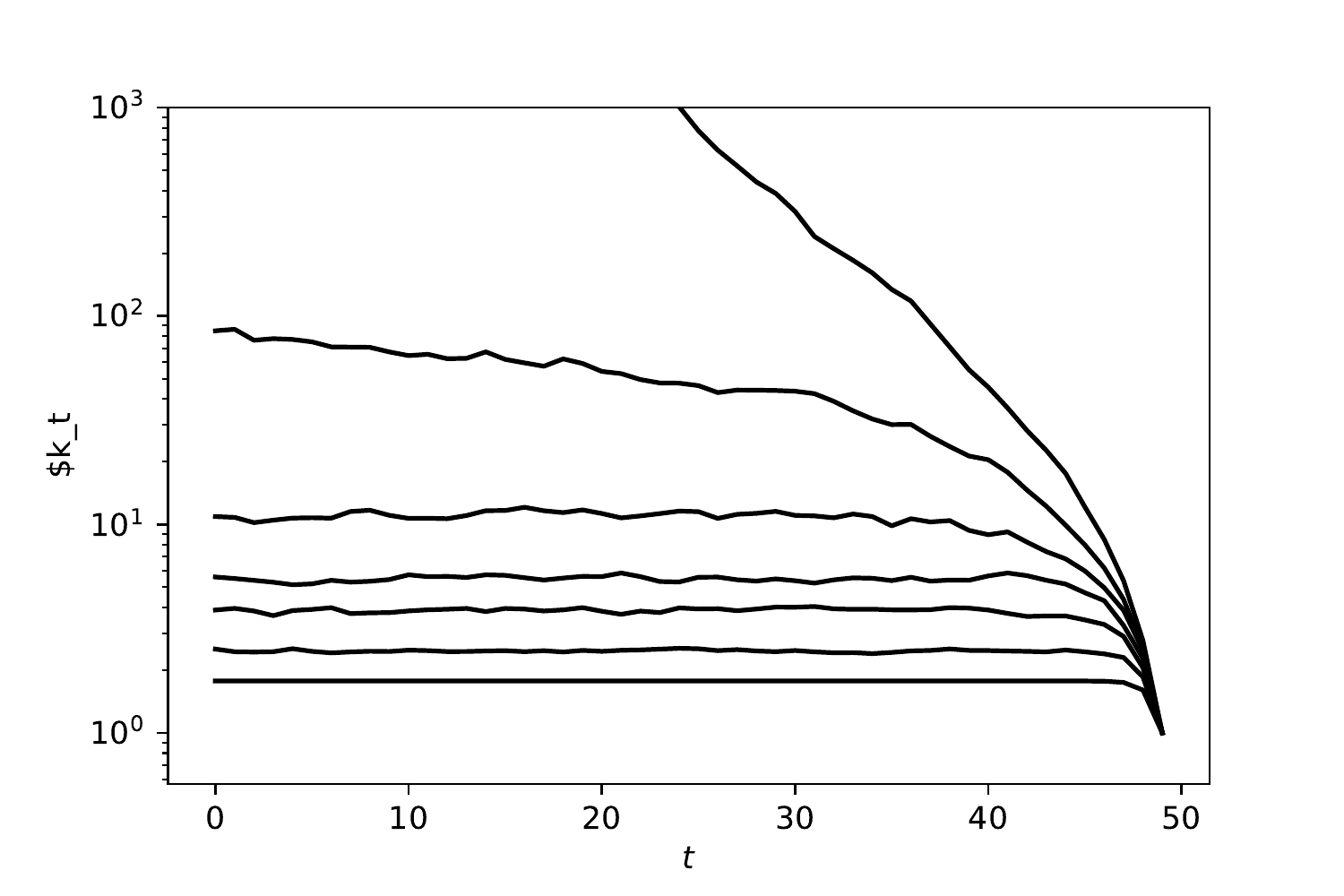}
	    \caption{}
	    \label{fig:random2}
    \end{subfigure}
    \begin{subfigure}{.4\textwidth}
	    \includegraphics[width=\linewidth]{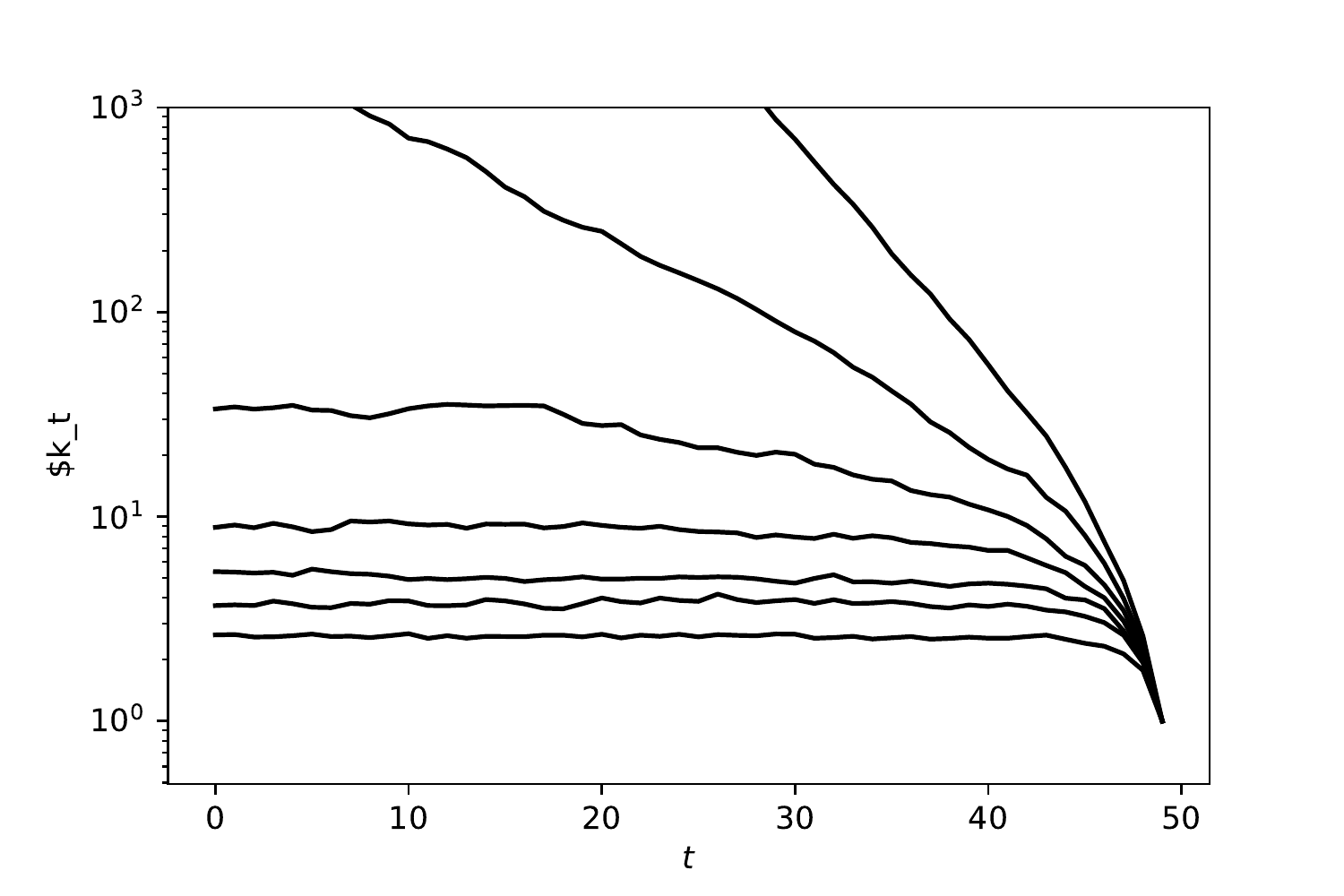}
	    \caption{}
	    \label{fig:random3}
    \end{subfigure}
    \caption{Random LQR example. (a) $\Sigma_{bb}=0$ and varying $\Sigma_{aa}$ (b) $\Sigma_{aa}=0$ and varying $\Sigma_{bb}$ (c) $\Sigma_{bb}=0.64$ and varying $\Sigma_{aa}$.}
\end{figure}

For all of the following examples, we let $\bar a = 1.1$, $\bar b = 1.0$, $T=50$, $N=50$, and we plot the coefficient $k_t$ versus $t$ in three cases:
\begin{itemize}
	\item \emph{$a$ fixed and $b$ random.} We fix $\Sigma_{aa}=0$, and vary
	\[\Sigma_{bb}\in\{0,0.81,1.44,2.25,2.89,3.61,4.41,4.84,5.76\}.\]
	The results are displayed in Figure~\ref{fig:random1}.
	\item \emph{$a$ random and $b$ fixed.}
	We fix $\Sigma_{bb}=0$, and vary \[\Sigma_{aa}\in\{0,0.25,0.49,0.64,0.81,1.00,1.21\}.\]
	The results are displayed in Figure~\ref{fig:random2}.
	\item \emph{Both $a$ and $b$ random.}
	We fix $\Sigma_{bb}=0.64$, and vary \[\Sigma_{aa}\in\{0,0.16,0.25,0.36,0.49,0.64,.81\}.\]
	The results are displayed in Figure~\ref{fig:random3}.
\end{itemize}

Our results match those of~\cite{athans1977uncertainty}, modulo Monte Carlo error.
When the variance gets too large, the cost-to-go functions diverge to $+\infty$, as predicted by the uncertainty threshold principle, and numerically checkable by our implementation.

\subsection{Jump LQR}

Jump LQR is LQR with dynamics that jump (or switch) between modes.
In infinite-horizon jump LQR, the dynamics are
\BEAS
x_{t+1} &=& A^{s_t} x_t + B^{s_t} u_t + c^{s_t}, \quad t=0,1,\ldots, \\
s_{t+1} &=& \text{$i$ with probability $\Pi_{ij}$ if $s_t=j$}, \quad t=0,1,\ldots,.
\EEAS
As in LQR, we adopt a quadratic form stage cost for simplicity.
For this problem, there is an optimal policy that is affine in $x$ for each mode.
We can find this policy using the algorithms described in this paper.
When there is no switching ($\Pi=I$), the problem reduces to $K$ separate LQRs (one for each mode).
Depending on the exact problem data and distributions, the policies found with switching can be substantially different than the non-switching policies.

\paragraph{Numerical example.}
Consider the one-dimensional system, \ie, $m=n=1$,
with dynamics described by
\BEAS
x_{t+1} &=& 1.2x_t + 0.1u_t; \quad s=1,\\
x_{t+1} &=& 0.8x_t - 0.1u_t; \quad s=2,
\EEAS
and Markov chain switching probabilities given by
\[
\Pi = \begin{bmatrix} .8 & .2 \\ .2 & .8 \end{bmatrix}.
\]
Consider the following stage cost function
\[
g(x,u) = \frac{1}{2}x^2 + \frac{1}{2}u^2.
\]

We solved the undiscounted ($\gamma=1$) infinite-horizon problem with $20$ applications
of the Bellman operator.
The optimal policy for the case with switching is
\[
\phi^s(x)=
\begin{cases}
    -2.541x & s=1\\
    0.919x & s=2.
\end{cases}
\]
We also solved the problem with no switching, resulting in the optimal policy
\[
\phi_\text{ind}^s(x)=
\begin{cases}
    -3.844x & s=1\\
    0.207x & s=2.
\end{cases}
\]

The two policies are substantially different.
The policies, over $100$ time steps in the actual system starting at $x_0=(10)$ and $s_0=1$, achieve an expected
cost of $J(\phi) = 16.16$ and $J(\phi_\text{ind})=18.53$ respectively
(averaged over $100$ simulations with the same random seed).
The policy found taking into account the switching
outperforms the policy that ignores the switching.



\subsection{Multi-mission LQR}

LQR can be extended to handle randomly switching costs.
In multi-mission LQR, the mode $s_t$ corresponds to the current cost assigned to the controller.
The dynamics are the same as LQR, but the costs depend on the mode, or
\[
g_t^s(x,u) = \frac{1}{2}\begin{bmatrix}x \\ u \\ 1\end{bmatrix}^T
G_t^s \begin{bmatrix}x \\ u \\ 1\end{bmatrix}.
\]
The algorithms described in this paper can be used to solve the finite or
infinite-horizon stochastic control problems, and result in
an affine policy for each ``mission''.

\paragraph{Numerical example.}
We apply the example above to a tracking mission.
Let $p_t\in\reals^2$ denote the position and $v_t\in\reals^2$ denote the
velocity of a point mass in two dimensions.
The state is $x_t=(p_t,v_t)$ and
the force applied is $u_t\in\reals^2$.
Suppose the dynamics are described by
\[
x_{t+1} =
\begin{bmatrix} 1&0&0.05&0\\0&1&0&0.05\\0&0&0.98&0\\0&0&0&0.98\end{bmatrix}x_t
+
\begin{bmatrix}0&0\\0&0\\0.05&0\\0&0.05\end{bmatrix}u_t.
\]

\begin{figure}
    \centering
    \begin{subfigure}{.5\textwidth}
	    \includegraphics[width=\linewidth]{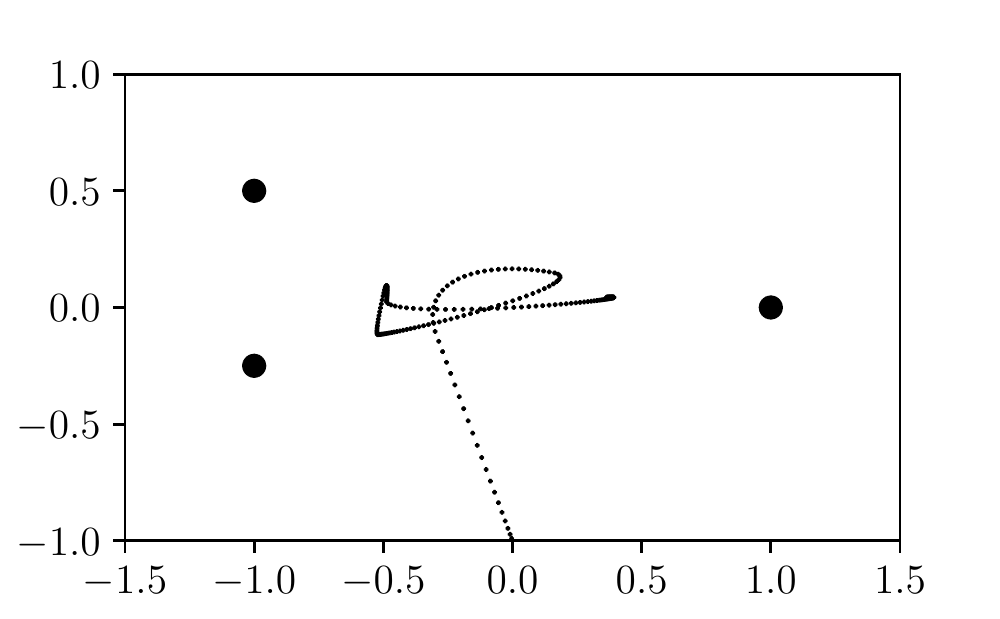}
	    \caption{average cost $=0.522$}
	    \label{fig:multi2}
    \end{subfigure}
    \begin{subfigure}{.5\textwidth}
	    \includegraphics[width=\linewidth]{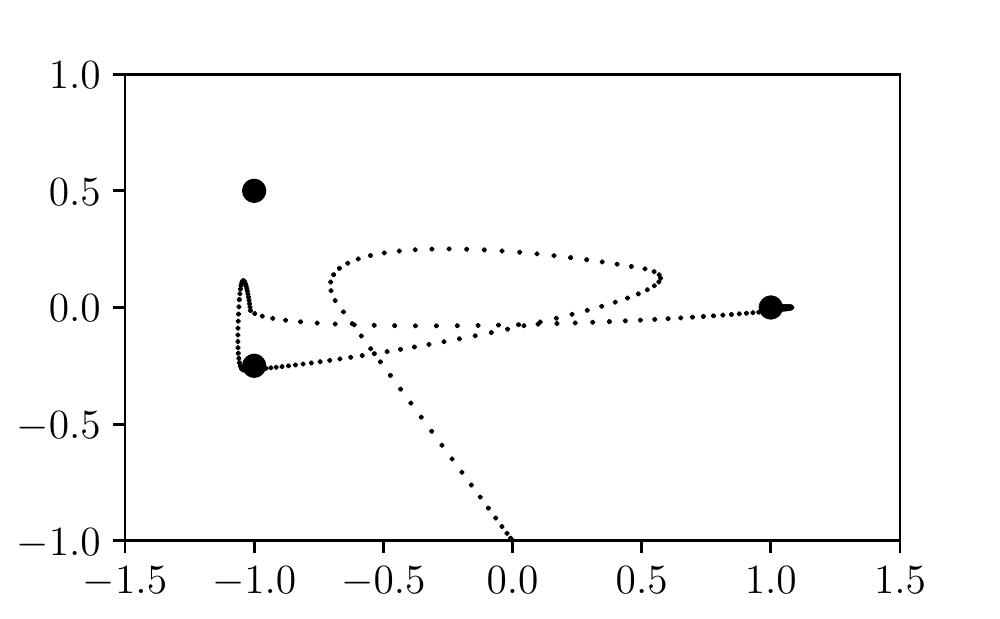}
	    \caption{average cost $=0.655$}
	    \label{fig:multi2}
    \end{subfigure}
    \caption{Multi-mission LQR example. (a) switching policy (b) non-switching policy.}
    \label{fig:multi}
\end{figure}

For each mission, the point mass's goal is to navigate to a target position $d_s\in\reals^2$
while minimizing control effort, corresponding to a stage cost
\[
g^s(x,u) = (1/2)\|p - d_s\|_2^2 + (\lambda/2)\|u\|_2^2,
\]
where $\lambda>0$.

Suppose we have three targets given by
\[
d_1 = (-1,.5) \quad d_2 = (-1,-.25) \quad d_3 = (1,0),
\]
and mission switching probabilities given by
\[
\Pi
=
\begin{bmatrix}
0.97 & 0.0075 & 0.015 \\
0.003 & 0.97 & 0.15 \\
0.027 & 0.0225 & 0.97
\end{bmatrix}.
\]
We solved the corresponding (infinite-horizon) stochastic control problem with $\lambda=0.1$, $T=50$, $N=1$, and $\gamma=1$.
A sample run of the optimal policy for the switching and no switching problems is shown in Figure~\ref{fig:multi}.
The policy that takes into account the switching does not go directly towards the targets, knowing that at any time the mission will switch and it will have to change course.





\subsection{Fault tolerant LQR}

We extend LQR to the case where control inputs (or actuators) randomly stop
affecting the system.
Suppose a system is deterministic and described by the linear dynamics
\[
x_{t+1} = Ax_t + Bu_t,
\]
and we would like to model input failures, \ie, if input $i$ has failed, $(u_t)_i$ has no effect on $x_{t+1}$.

We associate each mode of the system $s_i$ with an actuator configuration $a_i \subseteq \{1,\ldots,m\}$, which contains the indexes of the control inputs that have ``failed''. 
We assume that we know what actuator configuration we are in.
We can then associate each actuator configuration $a_i$ with a corresponding input matrix $B^i$, where
the $j$th column of $B^i$ is defined as
\[
B^i_j = \begin{cases} 0 & j \in a_i \\ b_j & \text{otherwise}\end{cases},
\]
where $b_j$ is the $j$th column of $B$.
We then can define a mode switching probability matrix $\Pi$, where $\Pi_{ij}$ is the probability that the system transitions from actuator configuration $a_j$ to actuator configuration $a_i$.
Our system's dynamics with actuator failures are described by
\BEAS
x_{t+1} &=& Ax_t + B^{s_t}x_t\\
s_{t+1} &=& i \; \text{with probability} \; \Pi_{ij} \; \text{if} \; s_t = j.
\EEAS

This is an extended quadratic control problem.
We also know that for each actuator configuration, the optimal policy is affine.
To implement this policy, one first identifies which actuators have failed, and then applies an affine transformation to the state to get the input.

This exact problem (with quadratic form cost) was first identified and solved by Birdwell
and Athans in 1977~\cite{birdwell1977reliability} (see also Birdwell's thesis~\cite{birdwell1978thesis}), where they derived what they call ``a set of highly coupled Riccati-like matrix difference equations''.
As expected, the cost-to-go functions that they derive are quadratic and the optimal policy is affine for each actuator configuration.
They used about a page of algebra to show this; we immediately know this.

\paragraph{Numerical example.}

We reproduce a slight modification of the example in~\cite{birdwell1977reliability}.
Here our dynamics are
\[
x_{t+1} = \begin{bmatrix}2.71828&0\\0&0.36788\end{bmatrix}x_t + \begin{bmatrix}1.71828&1.71828\\-0.63212&0.63212\end{bmatrix} u_t.
\]
The actuator configurations are $\emptyset$, $\{1\}$, and $\{2\}$, resulting in control matrices
\[
B_1 = \begin{bmatrix}1.71828&1.71828\\-0.63212&0.63212\end{bmatrix},
\quad B^2 = \begin{bmatrix}0&1.71828\\0&0.63212\end{bmatrix}, \quad, B^3 = \begin{bmatrix}1.71828&0\\-0.63212&0\end{bmatrix}.
\]
The mode switching probability matrix is
\[
\Pi = \begin{bmatrix}
0.943 & 0.069 & 0.026 \\
0.03  & 0.854 & 0.04 \\
0.027 & 0.077 & 0.934 \\
\end{bmatrix}.
\]
We use a stage cost $g(x,t) = \|x\|_2^2 + \|u\|_2^2$.
The optimal policies are linear, given by
\BEAS
\phi^1(x) &=& \begin{bmatrix}-0.737&0.135\\-0.74&-0.136\end{bmatrix}x \\
\phi^2(x) &=& \begin{bmatrix}0&0\\-1.455&-0.003\end{bmatrix}x \\
\phi^3(x) &=& \begin{bmatrix}-1.462&0.002\\0&0\end{bmatrix}x.
\EEAS
If an actuator does not affect the system, the optimal action sets this input to zero,
since one incurs cost for having its input not equal to zero.
(This is why the first column of $B^2$ is zero and the second column of $B^3$ is zero.)
If we ignore the fact that the dynamics switch when an actuator fails, 
we arrive at a suboptimal policy.

\subsection{Portfolio allocation with multiple regimes}



In this example we frame the problem of designing an optimal portfolio
allocation in a market that randomly switches between multiple regimes.
We borrow the notation from~\cite{boyd2016multi}.

\paragraph{Holdings.}
We work in a universe of $n$ financial assets.
We let $h_t\in\reals^n$ denote the (dollar-valued) holdings of our portfolio in each of
those $n$ assets at the beginning of time period $t$.
(We allow $(h_t)_i<0$, which indicates that we are short selling asset $i$.)
The total value of our portfolio at time $t$ is $v_t=\ones^T h_t$.
Our state is $h_t$.

\paragraph{Trading.}
At the beginning of each time period, we select a trade vector $u_t\in\reals^n$ that denotes the dollar value of the trades to be executed.
After making these trades, the investments are held constant until the next time period.
The post-trade portfolio is denoted
\[
h_t^+ = h_t + u_t, \quad t=0,\ldots,T-1,
\]
with a total value $v_t^+ = \ones^T h_t^+$.
We have a self-financing constraint, \ie, $v_t=v_t^+$, which can be expressed as
\[
\ones^Tu_t=0.
\]

\paragraph{Market state.}

We assume that the market
is in one of several (fully observable) market regimes.
We assign a mode $s_t\in\{1,\ldots,K\}$ to each regime.
Each regime corresponds to a different return distribution and transaction costs.
We use a Markov chain to model the market regime switching from time $t$ to $t+1$.

\paragraph{Investing dynamics.}
The post-trade portfolio is invested for one period.
Assuming the market is in the mode $s_t$, the portfolio at the next time period is given by
\[
h_{t+1} = (I+\diag(r_t^{s_t})) h_t^+,
\]
where $r_t^{s_t}\in\reals^n$ is a random vector of asset returns from time $t$ to time $t+1$ when the market is in the mode $s_t$.
The mean of the return vector for time $t$ in mode $s_t$ is denoted $\mu_t^{s_t} = \Expect[r_t^{s_t}]$ and its covariance is denoted
$\Sigma_t^{s_t} = \Expect[(r_t^{s_t}-\mu_t^{s_t})(r_t^{s_t}-\mu_t^{s_t})^T]$.

\paragraph{Transaction cost.}
The trading results in a transaction cost $\phi_{t,s}^\text{trade}:\reals^n \rightarrow\reals$ (in dollars), which we assume to be a (diagonal) quadratic form, \ie,
\[
\phi_{t,s}^\text{trade}(u_t) = u_t^T\diag(b_t^{s})u_t,
\]
where $(b_t^{s})_i\in\reals^+$ is the coefficient of the quadratic transaction cost for asset $i$ during time period $t$ when the market is in the mode $s$.

\paragraph{Returns.}
The portoflio return from period $t$ to $t+1$ when the market is in mode $s$ is given by
\BEAS
R_t^s &=&v_{t+1}-v_t - \phi_{t,s}^\text{trade}(u_t)\\
&=& \ones^T (h_{t+1}-h_t) - \phi_{t,s}^\text{trade}(u_t)\\
&=& \ones^T \left((I+\diag(r_t^{s}))(h_t+u_t) - h_t\right) - \phi_{t,s}^\text{trade}(u_t)\\
&=& (r_t^{s})^T x_t+ \ones^T u_t+(r_t^{s})^Tu_t - \phi_{t,s}^\text{trade}(u_t)\\
&=& (r_t^{s})^T (x_t+u_t) - \phi_{t,s}^\text{trade}(u_t).
\EEAS
The expected return is given by
\[
\Expect[R_t^s] = (\mu_t^s)^T (x_t + u_t) - \phi_{t,s}^\text{trade}(u_t).
\]
The variance of the return is given by
\[
\mathbf{Var}[R_t^s] = (x_t + u_t)^T \Sigma_t^s (x_t + u_t).
\]

\paragraph{Stage cost function.}
Our goal will be to maximize a weighted combination of the mean and variance of the returns, while satisfying the self-financing condition.
This can be accomplished with the following (extended quadratic) stage cost function
\[
g_t^{s}(x,u) = \begin{cases}
-\Expect[R_t^s] + \phi_t^\text{trade}(u) + \gamma_t \mathbf{Var}[R_t^s] & \ones^Tu = 0 \\
+\infty & \text{otherwise}
\end{cases},
\]
for some parameter $\gamma_t > 0$ that trades off return and risk.

The above problem can be solved with the algorithms described in this paper by providing an oracle that provides samples of $r_t^s$ for all $t,s$ (we can estimate the mean and covariance from this), the transaction cost vector $b_t^s$ for all $t,s$, and the mode switching probability matrix $\Pi_t$ for all $t$.
The optimal trade vector when the market is in mode $s$ will be an affine function of the holdings, \ie,
\[
u_t = K_t^{s_t}x_t + k_t^{s_t}.
\]
It turns out that the optimal policies can be written in the following more interpretable form
\[
u_t = K_t^{s_t}(h_t - (h^\star)_t^{s_t}),
\]
where $(h^\star)_t^{s_t}=-(K_t^{s_t})^\dagger k_t^{s_t}$ is the \emph{desired} holdings vector in the regime $s_t$.
So, to select a trade vector, we calculate the difference between our current holdings and our desired holding, and then multiply that difference by a feedback gain matrix.

To the best of our knowledge, this application has not yet appeared in the literature.


\paragraph{Numerical example.}

\begin{figure}
    \begin{subfigure}{.5\textwidth}
        \centering
        \includegraphics[width=\linewidth]{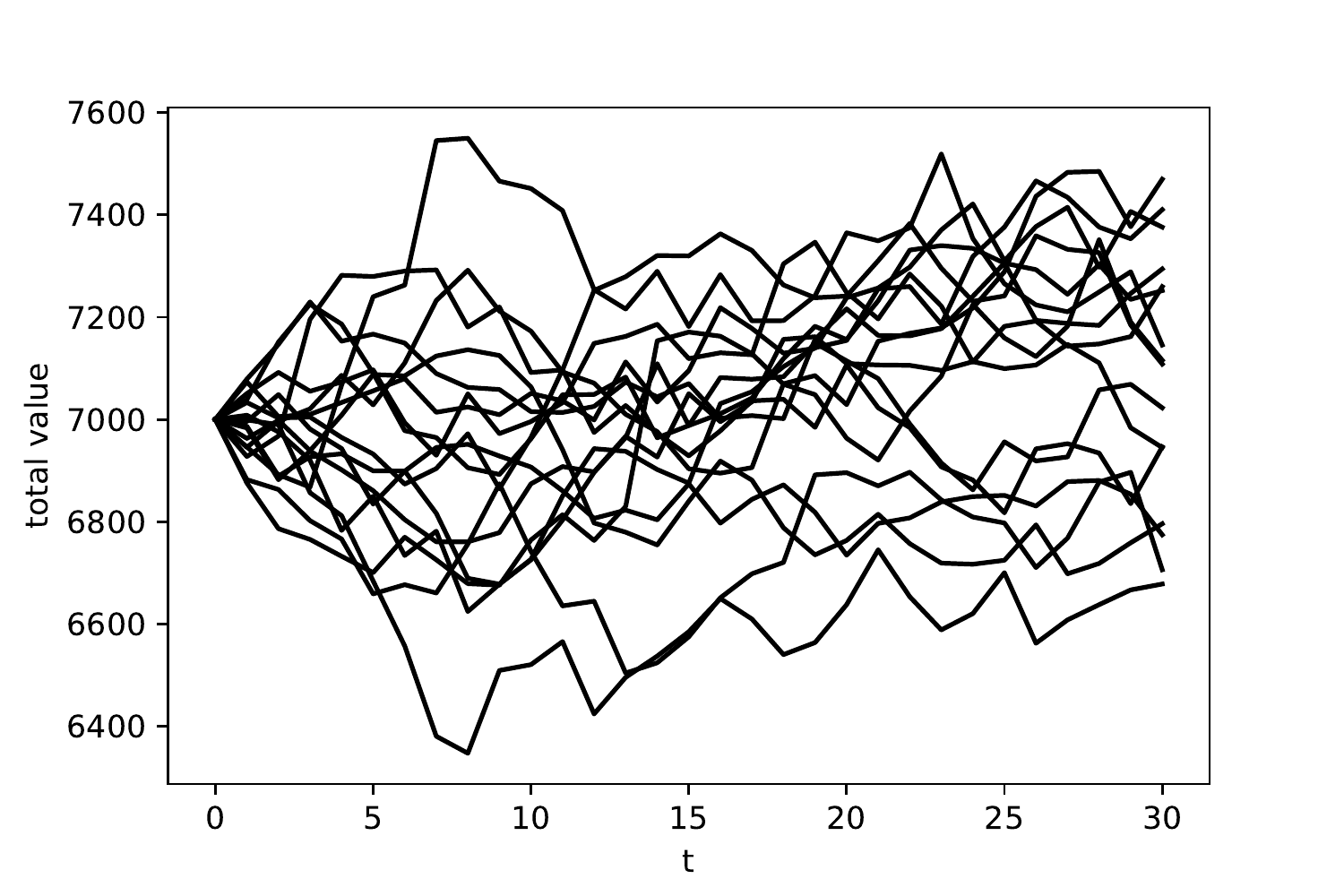}
        \caption{}
    \end{subfigure}%
    \begin{subfigure}{0.5\textwidth}
        \centering
        \includegraphics[width=\linewidth]{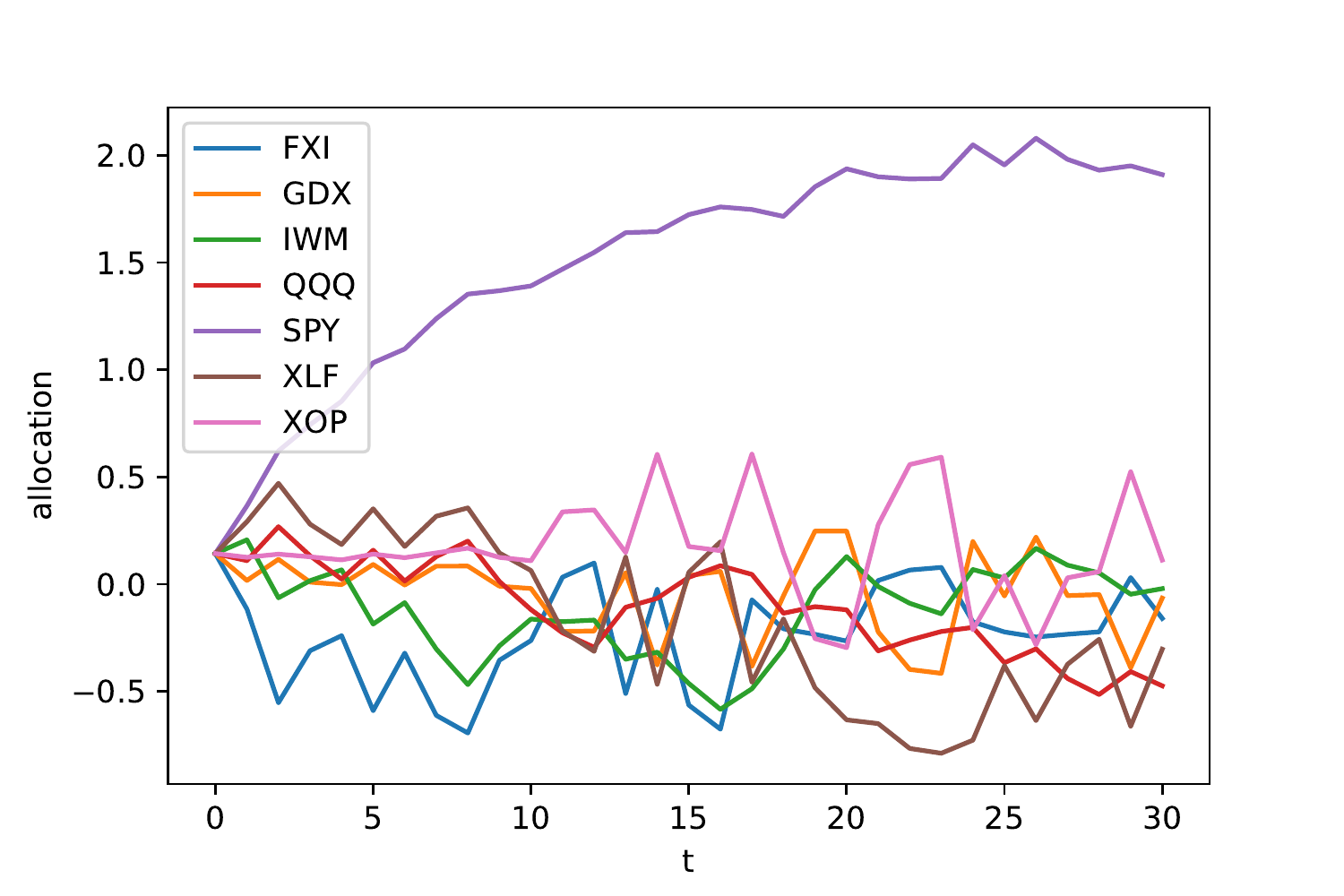}
        \caption{}
    \end{subfigure}
    \begin{subfigure}{.5\textwidth}
        \centering
        \includegraphics[width=\linewidth]{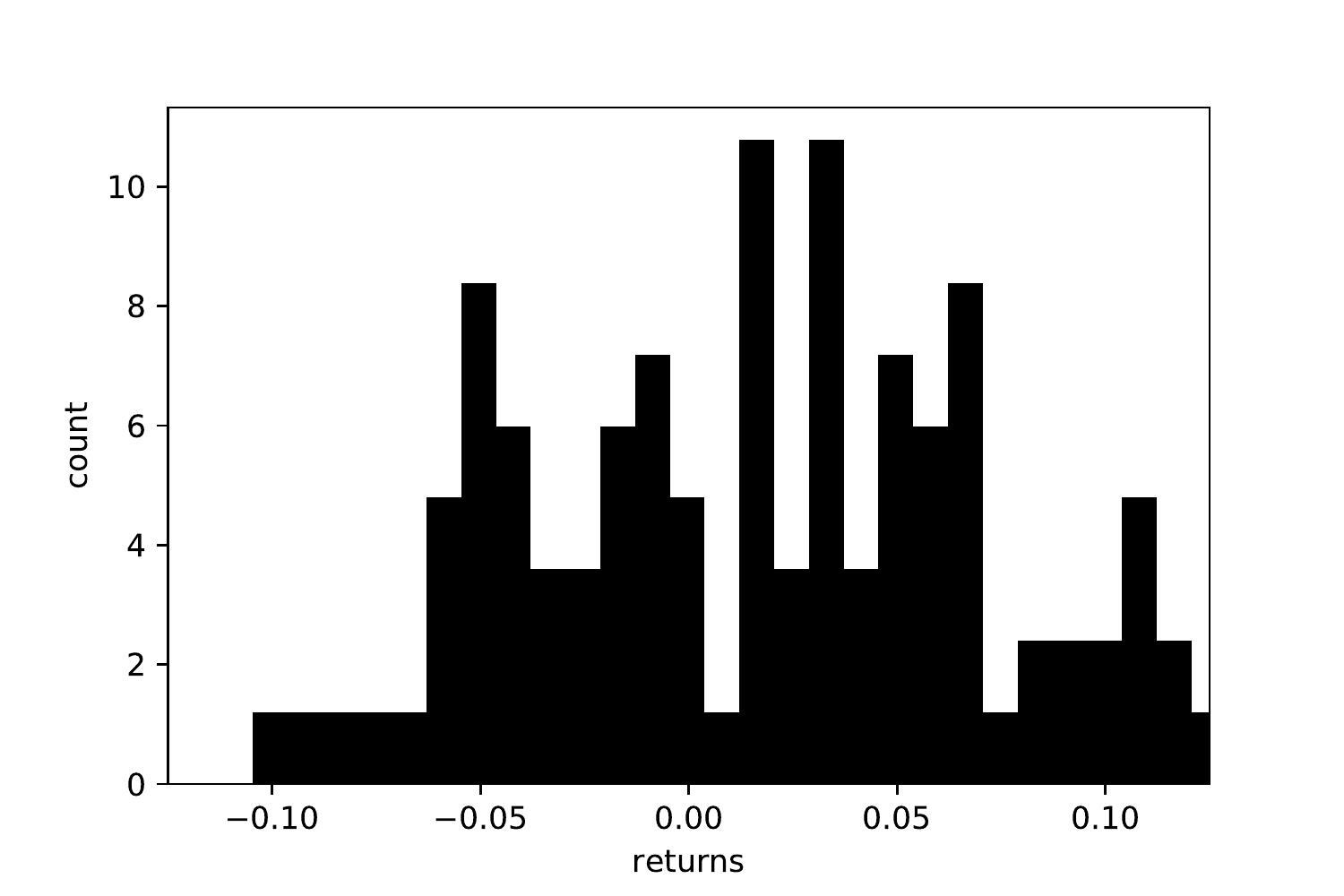}
        \caption{}
    \end{subfigure}%
    \begin{subfigure}{0.5\textwidth}
        \centering
        \includegraphics[width=\linewidth]{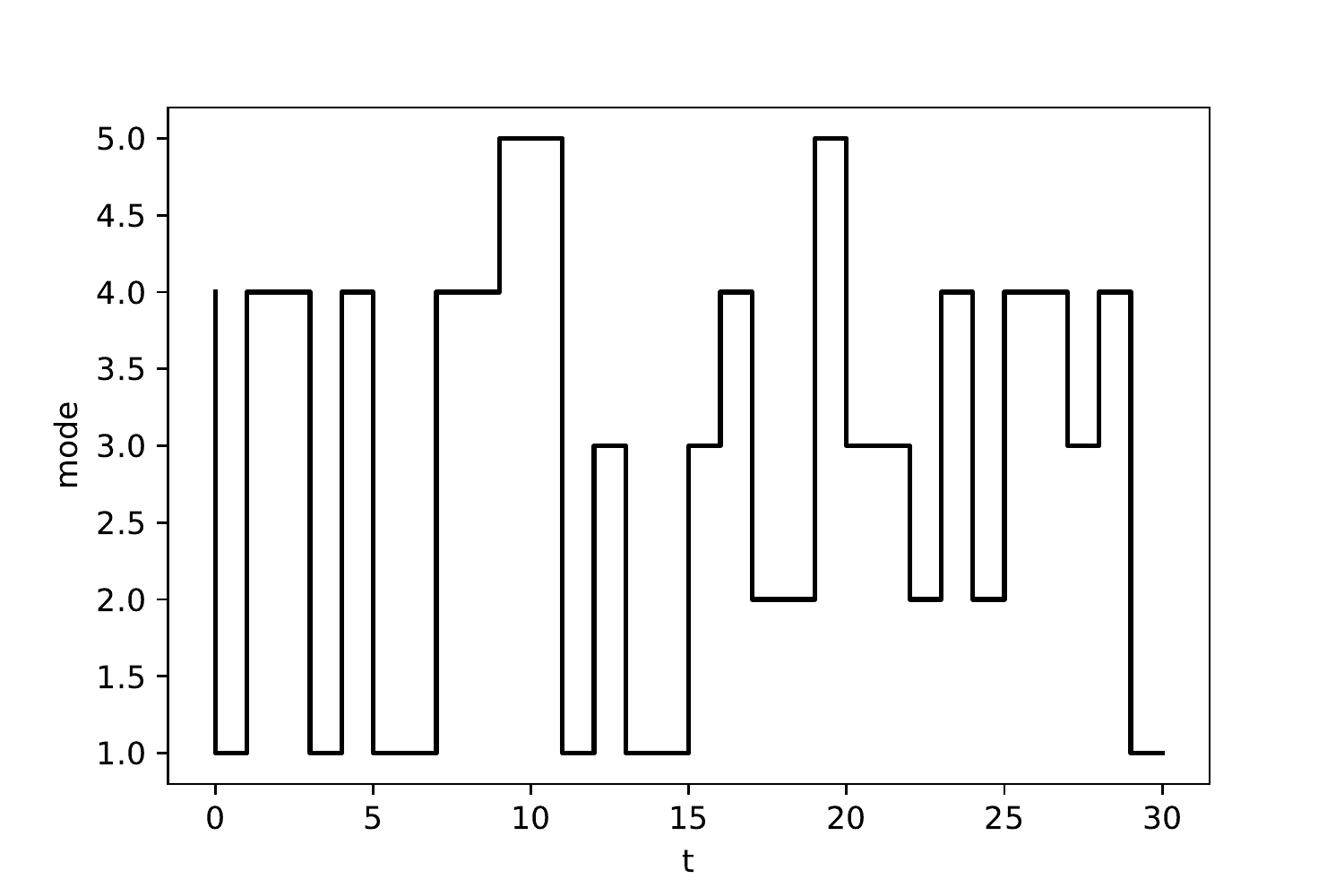}
        \caption{}
    \end{subfigure}
    \caption{(a) Total value of the portfolio versus time (in days) over fifteen simulations from the same initial condition. (b) Portfolio allocation ($h_t/v_t$) versus time (in days) for one simulation. (c) Histogram of returns over thirty days of trading over $1000$ simulations. The average return was $1.72\%$ with a standard deviation of $5.53\%$. (d) Market regime versus time for the same simulation as (b).}
    \label{fig:1}
\end{figure}

We gathered the daily returns from October 2013 to October 2018 of $n=m=6$ popular exchange traded funds (ETFs).

For the market regime, we used the daily rate of change of the CBOE Volatility Index (VIX), which approximates the market's expectation of 30-day volatility.
We gathered the daily opening price of the VIX and calculated its daily rate of change, which we refer to as dVIX.
We segmented dVIX into $K=5$ numerical ranges defined by the endpoints
\[
(-0.09 , -0.017, -0.003,  0.003,  0.03,  0.287),
\]
and define the range that it is in at time $t$ as the regime $s_t\in\{1,\ldots,K\}$.
We then calculated the empirical probabilities of switching between each market regime (values of dVIX), resulting in the following mode switching probability matrix
\[
\Pi = \begin{bmatrix}
0.159 & 0.123 & 0.146 & 0.189 & 0.282 \\
0.242 & 0.291 & 0.299 & 0.276 & 0.197 \\
0.108 & 0.215 & 0.201 & 0.156 & 0.155 \\
0.357 & 0.318 & 0.305 & 0.319 & 0.225 \\
0.134 & 0.054 & 0.049 & 0.06  & 0.141
\end{bmatrix}.
\]

For each regime, we gathered all of the days where the market was in that regime and fit a multivariate lognormal distribution to $\ones + r_t^s$ for the $5$ ETFs (see Appendix~\ref{sec:lognormal}).

Using these distributions and mode switching probabilities, we solved an instance of the portfolio allocation problem with a time horizon of $T=30$, $N=50$, $\gamma_t=\num{1e-1}$, $b=p_0 \cdot \num{1e-7}$ (where $p_0\in\reals^n$ is the price of the assets at the final day of the ETF data), and no final cost.
We then simulated the system several times (using different random seeds) starting in the initial state $h_0=(1000)\ones$ and $s_0=3$.
Figure~\ref{fig:1} displays various quantities over time from the simulations.

\paragraph{Incorporating constraints.}

\begin{figure}
    \centering
    \includegraphics[width=.5\linewidth]{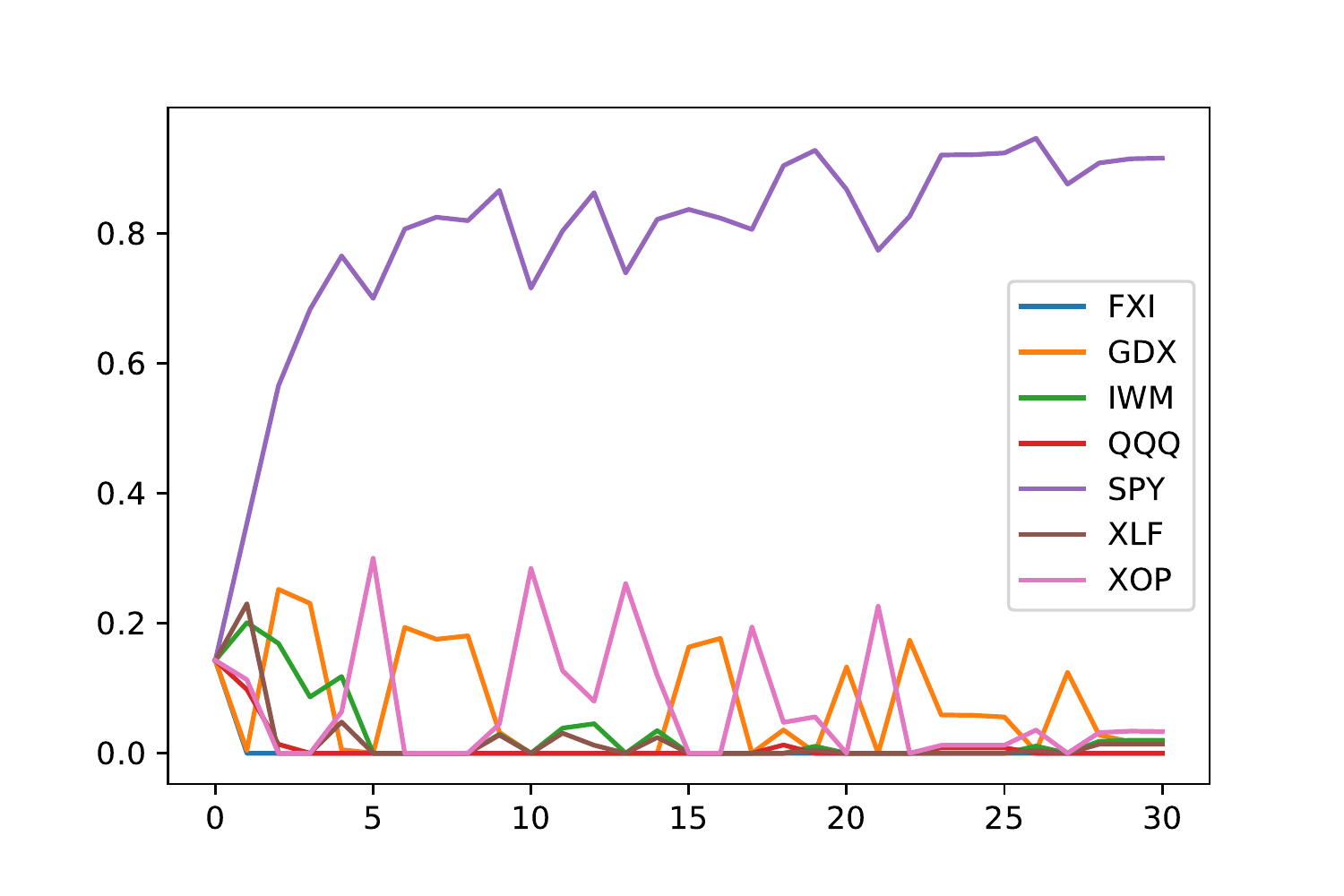}
    \caption{Asset allocation versus time for long-only portfolio. The average return over $12$ random seeds is $0.31\%$ with a standard deviation of $3.26\%$.}
    \label{fig:longonly}
\end{figure}

We can incorporate convex constraints on the trade vector using the technique described in \S\ref{sec:extandvars}.
For example, we can make the portfolio \emph{long-only}, meaning we constrain
\[
h_t+u_t \geq 0.
\]
For other useful convex constraints on $h_t$ and $u_t$, see~\cite{boyd2016multi}.
We cannot solve the stochastic problem exactly with these constraints, but we can approximate the solution by instead minimizing our approximate state-action cost-to-go function with the constraints on the trade vector (\ie, by solving a convex quadratic program).
To evaluate the policy, we must solve a convex quadratic program, which is slower than just a matrix-vector multiply.
(There are methods which can speed up evaluation of such policies, see, \eg,~\cite{wang2011fast}.)
In Figure~\ref{fig:longonly} we plot the asset allocation over time in the long-only case for the same random seed as the non-constrained example.


\subsection{Optimal execution}
\label{sec:optexec}

In this example, we frame the problem of designing an optimal liquidation scheme
for a financial asset, while minimizing transaction costs
and exposure, as a stochastic control problem.

We let $q_t\in\reals^n$ denote the quantity (in shares) of $n$ assets at time $t$ that we own and would like to liquidate (get rid of) by time $T$.
Let $p_t\in\reals^n$ denote the price (in dollars) of each of these assets
at the beginning of time $t$.
Our state is $x_t=(q_t,p_t)$.
At time $t$, we sell $u_t\in\reals^n$ shares of
the asset ($(u_t)_i<0$ corresponds to buying asset $i$).
After making the trade, at time $t+1$, we own
\[
q_{t+1} = q_t - u_t
\]
shares.
We assume that our trade has linear market impact, resulting in the post-trade price
\[
p_t^+ = p_t - G^\text{imp}_t u_t,
\]
where $G^\text{imp}_t\in\reals^{n \times n}$ is the linear market impact matrix for time $t$.
The price is then updated from time $t$ to time $t+1$ according to
\[
p_{t+1} = A_t p_t^+,
\]
where $A_t\in\reals^{n \times n}$ is a (random) multiplicative total return matrix for time period $t$.
A simple model for $A_t$ is that it is diagonal.
The full dynamics are then
\BEAS
q_{t+1} &=& q_t-u_t\\
p_{t+1} &=& A_t(p_t-G^\text{imp}_tu_t).
\EEAS

At time $t$, we receive 
\[
u_t^T p_t - u_t^T \diag(\gamma^\text{tr}_t) u_t
\]
dollars from our transaction, where $\gamma^\text{tr}_t \in \reals^n$, which is element-wise positive,
is the (quadratic) transaction cost parameter for each of the assets at time $t$.
We would also like to minimize our exposure to the assets,
which we define as $q_t^T \diag(\gamma^\text{ex}_t) q_t$,
where $\gamma^\text{ex}_t\in\reals+^n$ is the element-wise positive 
exposure cost parameter at time $t$.
(This makes sense, because $q_t=0$ corresponds to zero exposure to the market.)

Our goal is to maximize profit while minimizing exposure, resulting in a (quadratic) stage cost of
\[
g_t(x,u,w) = -u^Tp + u^T\diag(\gamma^\text{tr}_t) u + q^T \diag(\gamma^\text{ex}_t) q
\]
with the constraint that $q_T=0$.
We can use the algorithms described in this paper to solve the optimal
execution problem if we are given the distributions for $r_t$
and the values of $G^\text{im}_t$, $\gamma^\text{tr}_t$, and $\gamma^\text{ex}_t$.

This problem was first identified and solved by Bertsimas and Lo in 1998~\cite{bertsimas1998optimal} for a single asset and later extended by the same authors to a (multiple asset) portfolio~\cite{bertsimas1999optimal} where the assets have correlated returns.

\paragraph{Numerical examples.}

We focus on the case $n=1$, \ie, we only hold one asset.
(We can easily solve the case where $n>1$, but it is harder to visualize.)
We use an i.i.d.\ log-normal distribution for the returns, or
\[
r_t\sim\text{Lognormal}(\mu,\sigma^2),
\]
for some $\mu\in\reals$ and $\sigma>0$.
We use a time horizon $T=7$, and coefficients that do not vary with time.

The simplest non-trivial example that we can construct is where 
\[\gamma^\text{tr} >0, \quad G^\text{imp}=\gamma^\text{ex}=\mu=\sigma=0.\]
Here the asset stays the same price, there is no exposure cost/market impact, and there is a quadratic transaction cost.
There is a simple analytical solution; the optimal policy is to sell a fixed fraction $u_t=q_0/T$ of the asset every time period, achieving
a total cost of $-q_0p_0 + \frac{1}{2}\gamma^\text{tr} q_0^2/T$.
With our implementation, we find the policy
\BEAS
\phi_0(q,p) &=& (0.14)q\\
\phi_1(q,p) &=& (0.17)q\\
\phi_2(q,p) &=& (0.2)q\\
\phi_3(q,p) &=& (0.25)q\\
\phi_4(q,p) &=& (0.33)q\\
\phi_5(q,p) &=& (0.5)q\\
\phi_6(q,p) &=& (1)q,
\EEAS
and the first cost-to-go function is
\[
V_0(q,p) = \frac{1}{2}\begin{bmatrix}q\\p\\1\end{bmatrix}^T\begin{bmatrix}0.14 & -1 \\ -1 & 0\end{bmatrix}\begin{bmatrix}q\\p\\1\end{bmatrix}.
\]
This matches what we expected.

\begin{figure}
    \begin{minipage}{.5\textwidth}
        \centering
        \includegraphics[width=\linewidth]{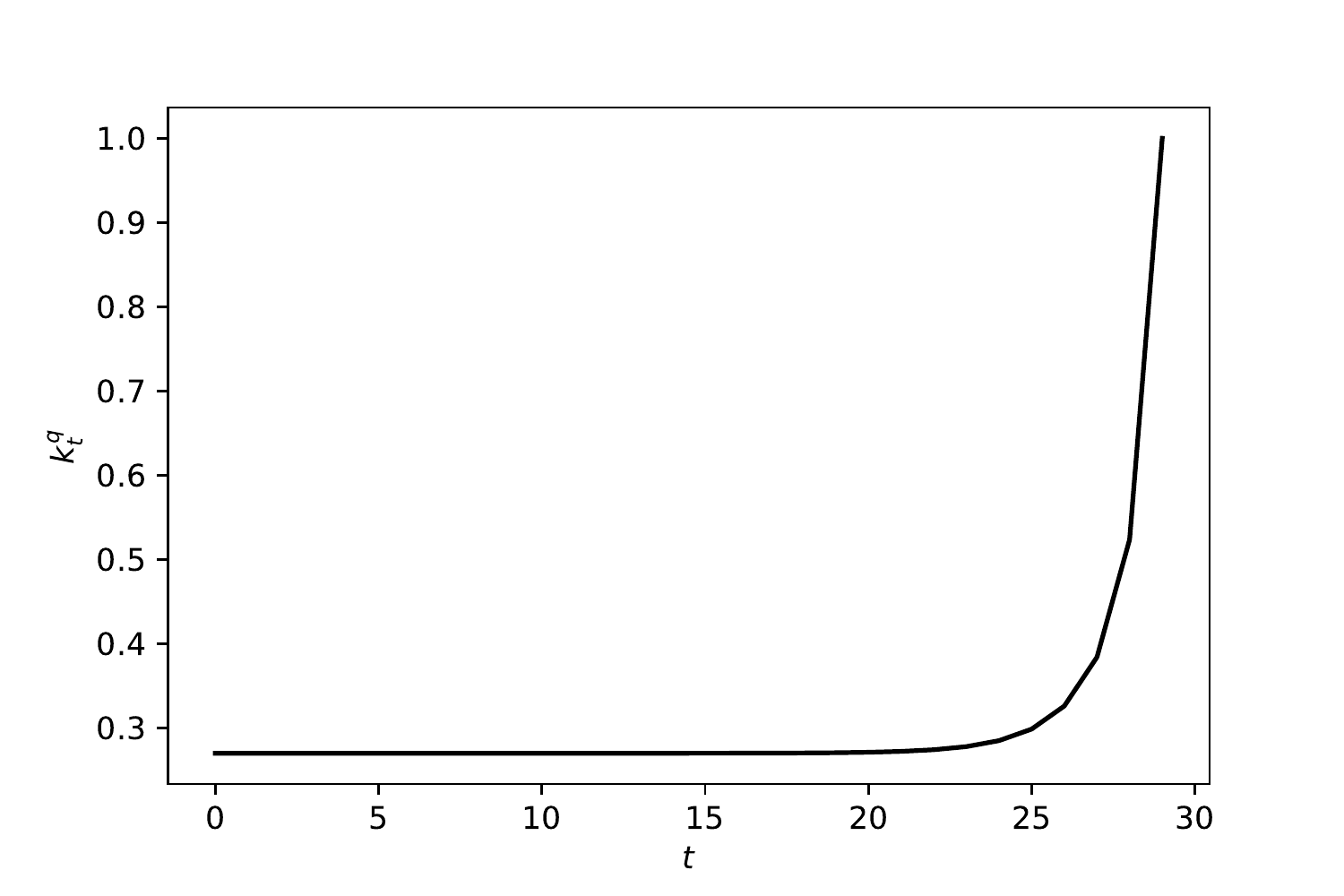}
        \caption{Coefficient for quantity.}
        \label{fig:ex1}
    \end{minipage}%
    \begin{minipage}{0.5\textwidth}
        \centering
        \includegraphics[width=\linewidth]{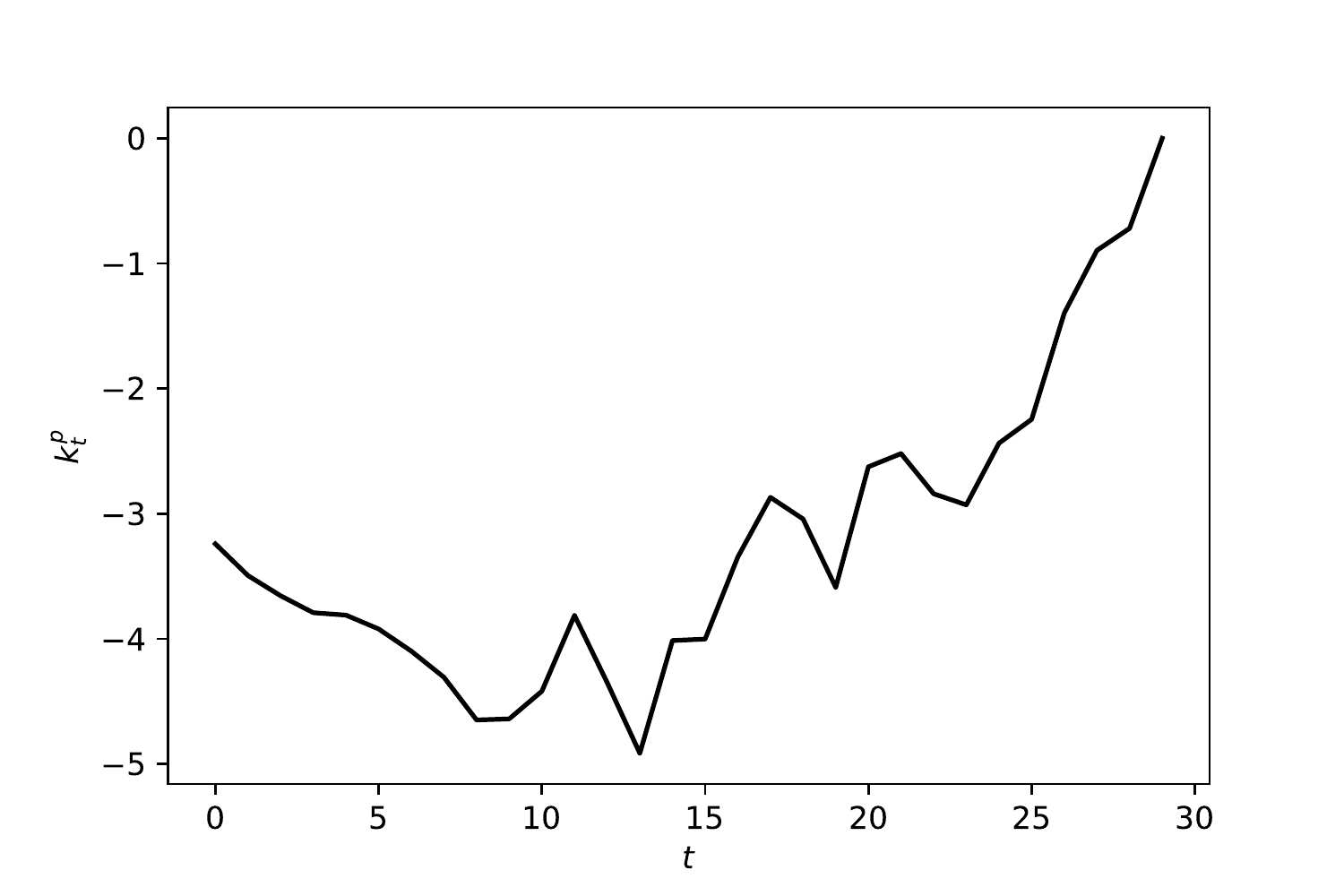}
        \caption{Coefficient for price.}
        \label{fig:ex2}
    \end{minipage}
\end{figure}

We gathered the daily returns of AAPL from January $2018$ to October $2018$, and fit a lognormal distribution to them.
We then solved the optimal execution problem over a time horizon $T=30$ (and $N=100$) with the parameters
\[
\gamma^\text{tr} =\num{1e-3}, \quad G^\text{imp}=\num{1e-4}, \quad \gamma^\text{ex}=\num{1e-4}, \quad \mu=0.0013, \quad \sigma=0.0153.
\]
The optimal policy has the form
\[
\phi_t(q,p) = k_t^q q + k_t^p p.
\]
Figures~\ref{fig:ex1} and~\ref{fig:ex2} show the coefficients of the policy versus time.
Interestingly, the higher the price of AAPL, the more stock one should buy (and this decreases linearly as time increases).

\begin{figure}
    \centering
    \includegraphics[width=.5\linewidth]{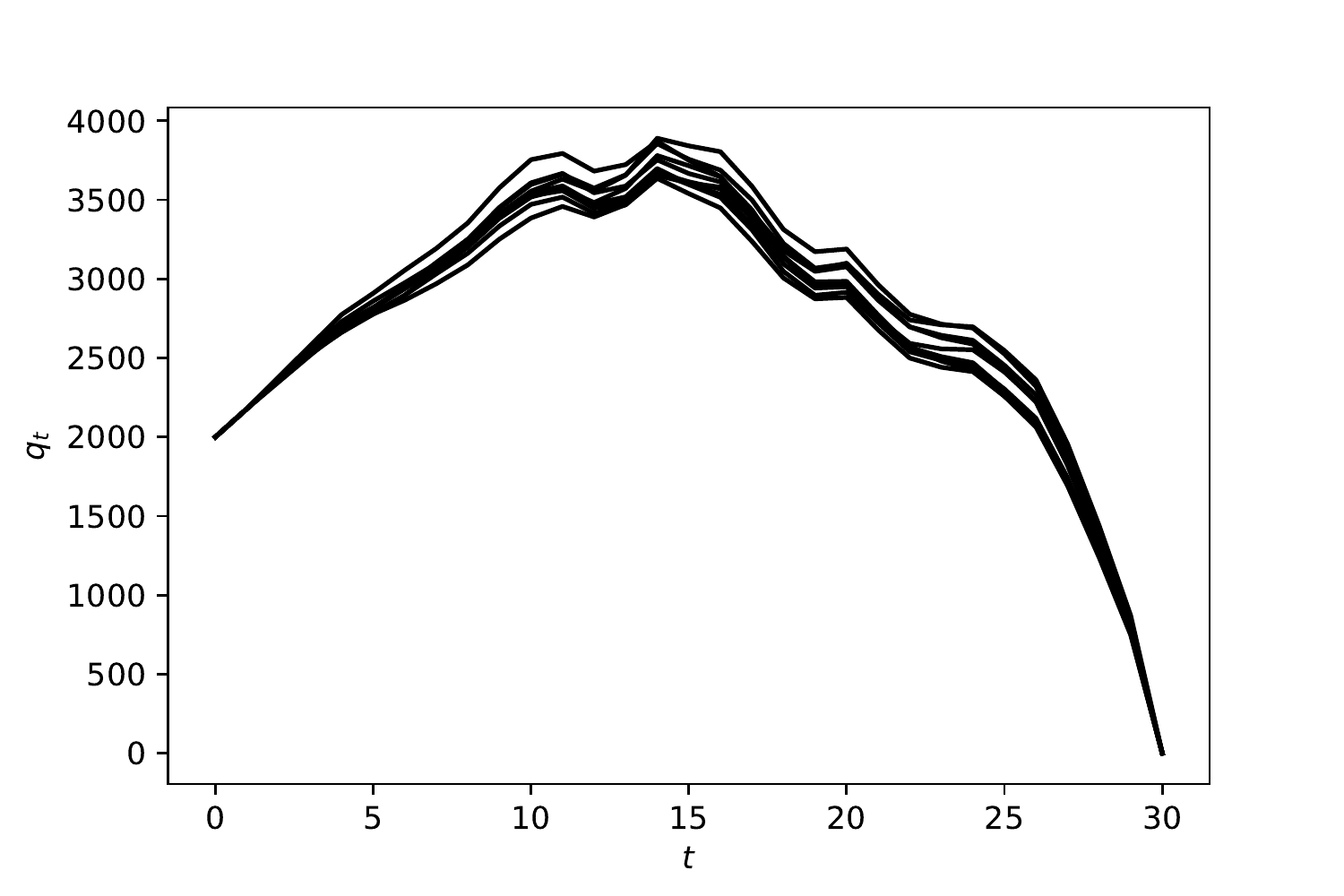}
    \caption{$q_t$ versus $t$ for $10$ random seeds.}
    \label{fig:ex3}
\end{figure}

We simulated the AAPL example over $10$ random seeds, with the same coefficients as above, in the initial state $(2000,222)$.
In other words, we wish to liquidate $2000$ shares of AAPL starting at $\$222$ over $30$ trading days (about six weeks).
In Figure~\ref{fig:ex3}, we show the quantity of AAPL stock held by the optimal policy over time.
At the beginning of the month, we buy some shares of AAPL, since it is likely to increase, then after about two weeks, we slowly begin to sell.
Taking into account transaction costs, this policy returns the investor $\$475\,513 \pm \$41\,614$, where the initial value of the investment was $\$444\,000$.

\subsection{Optimal execution with a random horizon}

We can adapt the above optimal execution example to make the time that we have to liquidate random.
At each time step $t$, we assume that there is a probability $p_t\in[0,1]$ that we will need to liquidate the asset at time $t+1$.
We add a mode $s_t\in\{1,2\}$; $s_t=2$ means we are forced to liquidate the asset, which we can enforce by adding the constraint $q_t=0$ to the stage cost when $s_t=2$.
The switching probability matrix is then
\[
\Pi_t = \begin{bmatrix}1-p_t & 0\\p_t & 1\end{bmatrix}.
\]
To the best of our knowledge, this application has not yet appeared in the literature.

\paragraph{Numerical example.}

\begin{figure}
    \begin{minipage}{.5\textwidth}
        \centering
        \includegraphics[width=\linewidth]{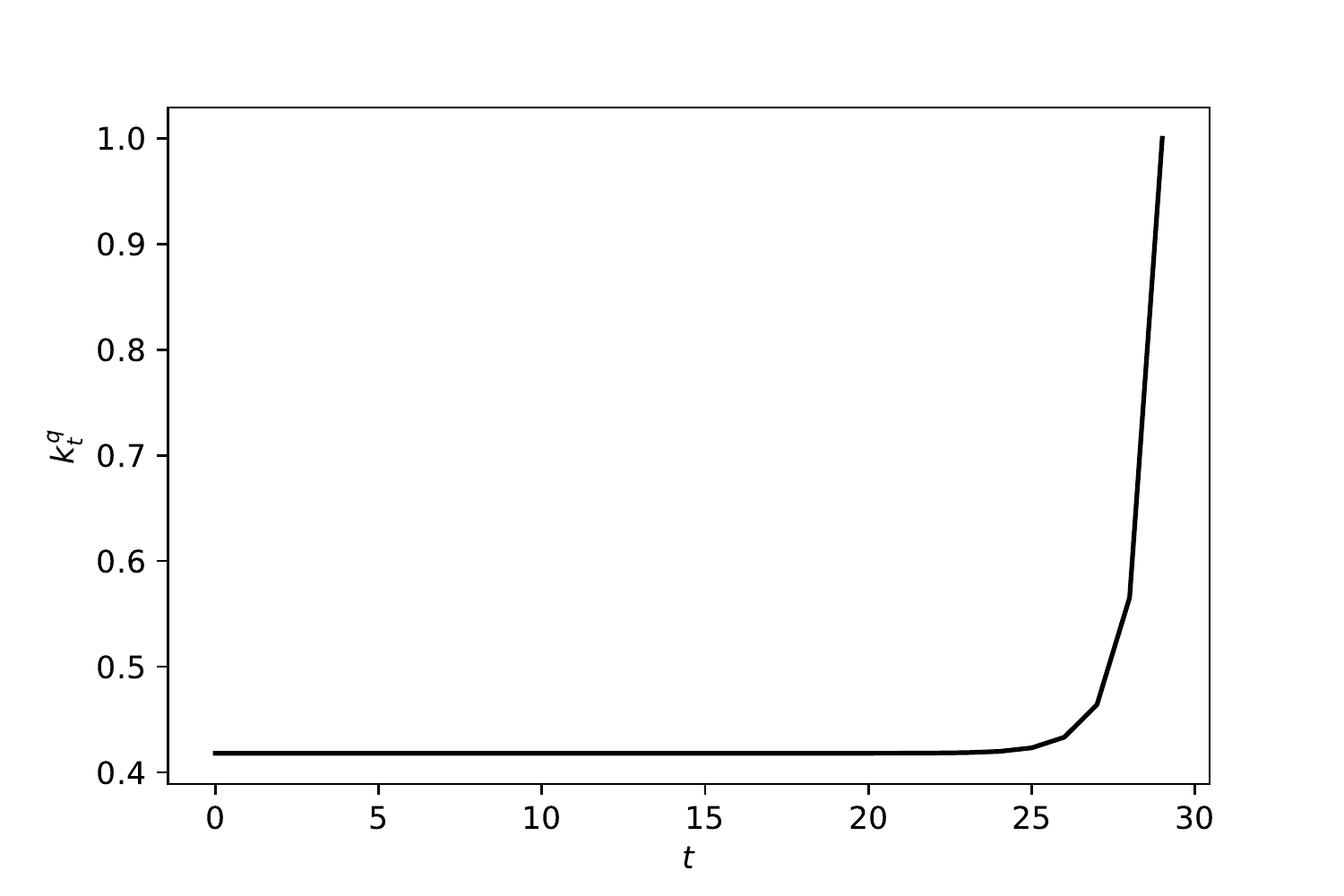}
        \caption{Coefficient for quantity.}
        \label{fig:exrand1}
    \end{minipage}%
    \begin{minipage}{0.5\textwidth}
        \centering
        \includegraphics[width=\linewidth]{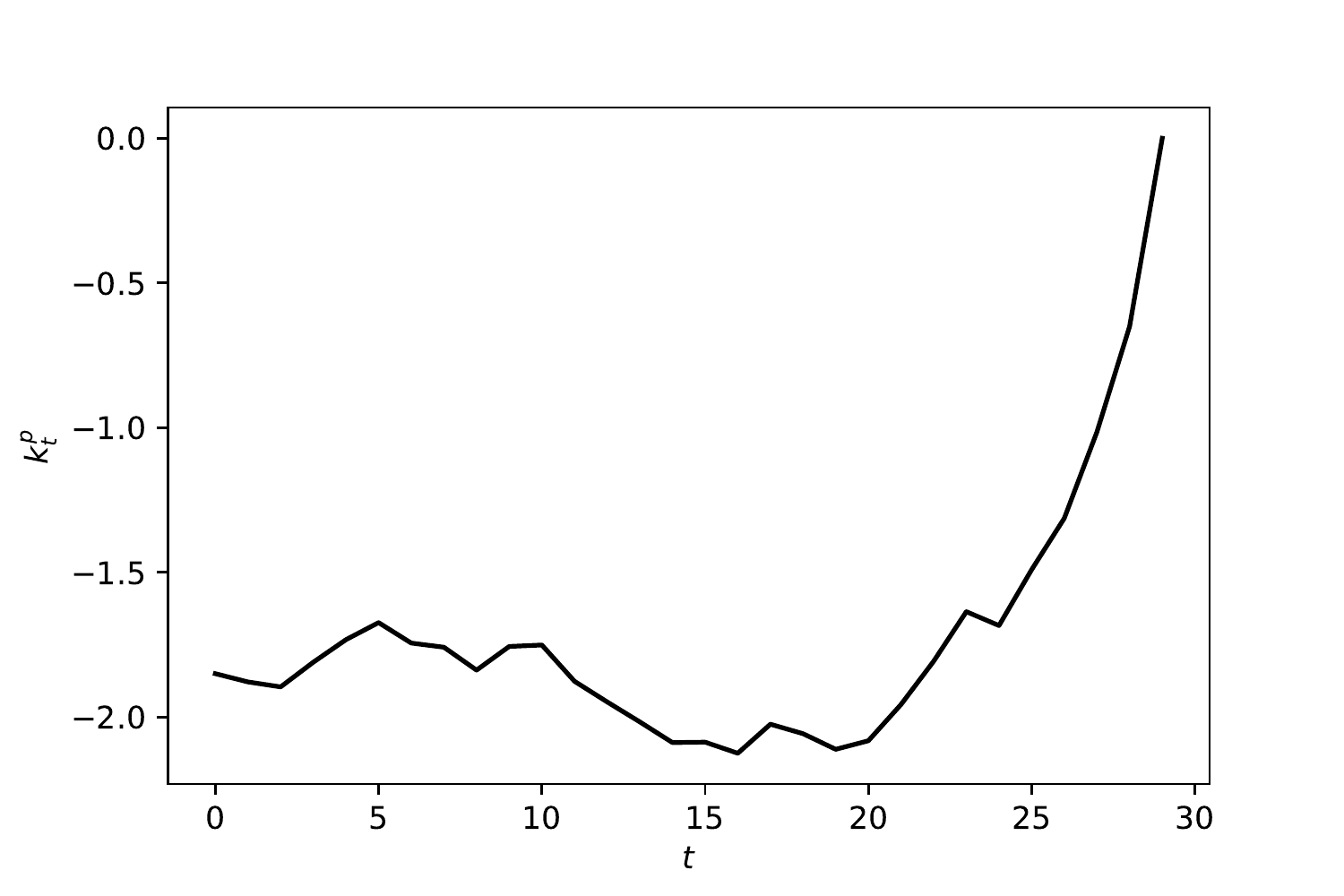}
        \caption{Coefficient for price.}
        \label{fig:exrand2}
    \end{minipage}
\end{figure}

\begin{figure}
    \centering
    \includegraphics[width=.4\linewidth]{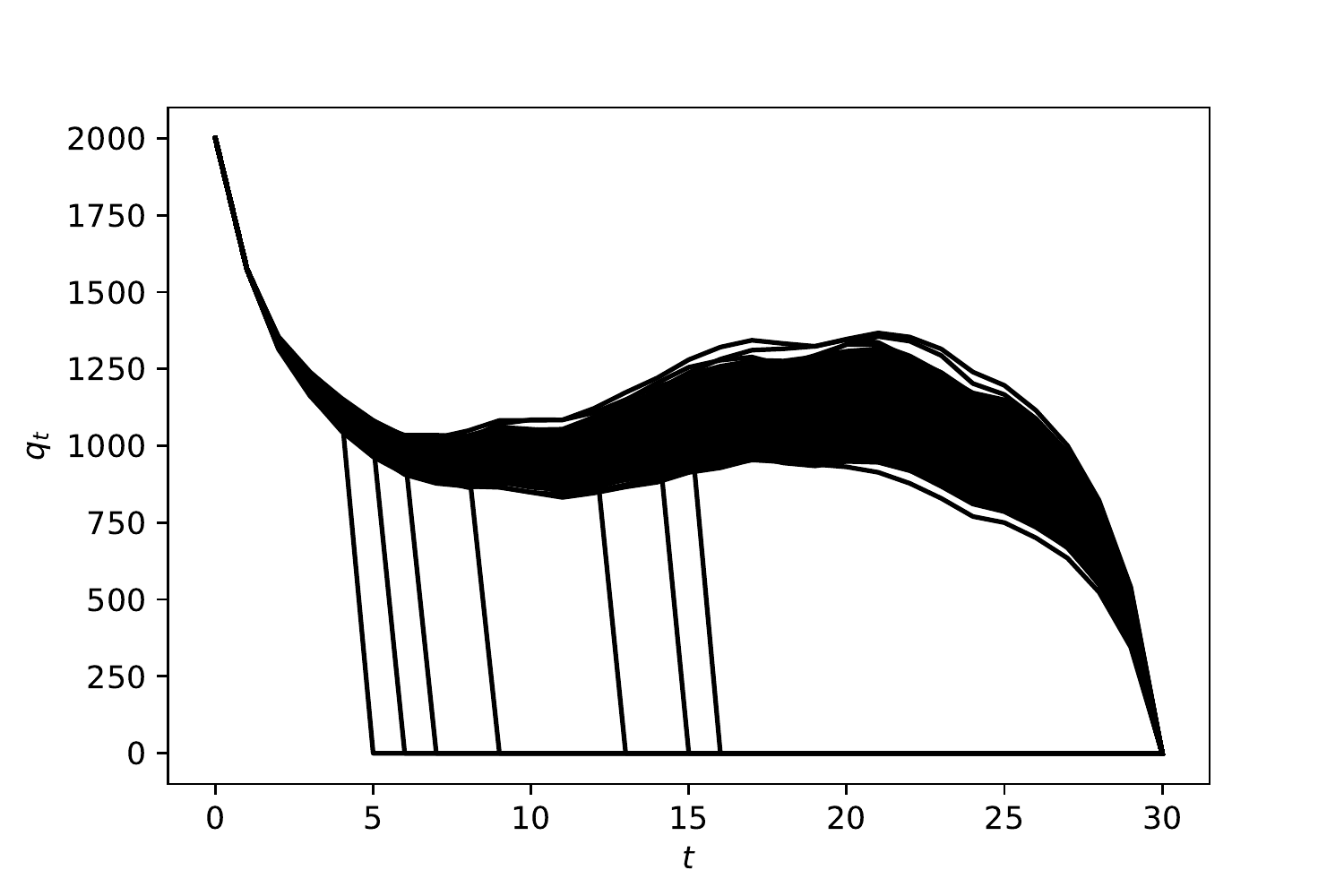}
    \caption{$q_t$ versus $t$ for $1000$ random seeds.}
    \label{fig:exrand3}
\end{figure}

We adapt the AAPL example in \S\ref{sec:optexec} by
adding a $.01\%$ probability that at each day, the next day we will have to fully liquidate the asset.
We produce the same figures as in the AAPL example in Figures~\ref{fig:exrand1}, \ref{fig:exrand2}, and \ref{fig:exrand3}.
The policy is much more cautious, because having to liquidate a large amount of stock in one day is very costly.

\subsection{Optimal retirement}

The goal in retirement planning is to devise an investment allocation
and consumption schedule for the rest of ones life in order to
maximize personal utility.
In this section we frame retirement planning as an extended quadratic control problem,
with the state being the investor's wealth,
the input being the allocation over various financial assets (and an amount to consume),
the cost being negative utility,
and the mode corresponding to whether the investor is alive or deceased.

We let the time period $t$ represent a (calendar) year.
At the beginning of year $t$, the mode $s_t$ corresponds to whether
the investor is alive ($s_t=1$)
or deceased ($s_t=2$).
If the investor is alive at the beginning of year $t$, either they pass away ($s_{t+1}=2$) with probability $p_t$, or they continue to live ($s_{t+1}=1$) with probability $1-p_t$.
The deceased mode is absorbing, \ie, if the investor is deceased at year $t$, they stay deceased at year $t+1$.
The mode dynamics is given by the mode switching probability matrix
\[
\Pi_t
=
\begin{bmatrix}
1-p_t & 0\\
p_t & 1\\
\end{bmatrix}.
\]

The investor's wealth $W_t\in\reals$ is the wealth (in dollars) of
the investor at the beginning of year $t$.
We operate in a universe of $m$ financial assets that the investor may choose to invest in.
At the beginning of year $t$, the investor allocates their wealth across
$n$ financial assets by specifying a trade vector $u_t\in\reals^m$, where $(u_t)_i$ is the amount (in dollars) bought of asset $i$ ($(u_t)_i<0$ corresponds to shorting the asset).
The amount that the investor does not invest,
\[
C_t = W_t - \ones^T u_t,
\]
is consumed if the investor is alive and bequeathed (\ie, left to
beneficiaries by a will) if the investor is deceased.
If the investor consumes $C_t$ during year $t$, they receive $
U_t(C_t)
$
utility, where $U_t:\reals\rightarrow\reals$ is a concave quadratic utility function for
consumption.
If the investor bequeaths $C_t$, they receive
$
B(C_t)
$
utility, where $B:\reals\rightarrow\reals$ is a concave quadratic utility function for
bequeathing.
(We enforce that all of the investor's wealth is bequeathed when they die via the constraint that $u=0$ when $s=2$.)

The investor's wealth at year $t+1$ is
\[
W_{t+1} = \begin{cases}
r_t^T u_t & s=1\\
0 & s=2
\end{cases},
\]
where $r_t\in\reals^m$ is a random total return vector for the financial assets over year $t$,
with $\Expect[r_t] = \mu_t$ and $\text{Cov}[r_t] = \Sigma_t$.
The variance of the investor's wealth at year $t+1$, assuming they are alive, is
\[
\text{Var}(W_{t+1}) = u_t^T \Sigma_t u_t.
\]

Our goal is to maximize utility, while minimizing risk, resulting in the stage cost functions
\BEAS
g_t^1(W,u) &=& -U_t(C) + \gamma u^T \Sigma_t u \\
g_t^2(W,u) &=& -B(C) + \begin{cases}0&u=0\\+\infty&\text{otherwise}\end{cases}.
\EEAS
We ignore transaction costs.

It is worth noting that this problem is small enough to be discretized and solved exactly (with up to $5$ or $6$ assets) with any dynamics and cost.
So this example is just an illustration, and it is only sensible to use extended quadratic control when there are many assets.

\paragraph{Numerical example.}

We gathered inflation-adjusted yearly returns for $m=3$ assets: the S\&P 500, 3-month treasury bills, and 10-year treasury bonds over the past $80$ years ($1938-2018$)~\cite{damodaran2018returns}.
We fit a multivariate lognormal distribution to the returns.

The mean and covariance of $r_t$ are
\[
\Expect[r_t] = \begin{bmatrix}
1.0854 \\ 1.0005 \\ 1.0169
\end{bmatrix}, \quad
\mathbf{Cov}[r_t] = \begin{bmatrix}
0.0316 & 0.00077 & 0.00003 \\
0.00077 & 0.00142 & 0.00139 \\
0.00003 & 0.00139 & 0.00673
\end{bmatrix}.
\]

For mortality rates, we use the Social Security actuarial life table~\cite{bell2005life}, which gives the death probability (probability of dying in one year) for each age, averaged across the United States's population.

We use the utility functions
\[
U_t(C) = -\frac{1}{2}(0.2)C^2 + (20)C, \qquad
B(C) = -\frac{1}{2}(0.002)C^2 + (4)C.
\]
Here $U_t(0)=B(0)=0$, the maximum of $U_t$ is at $100$ where $U_t(100)=1000$, and
the maximum of $B$ is at $2000$ where $B(2000)=4000$.
The utility functions were designed so that the investor ideally consumes
\$$100$k every year, and bequeaths \$$2$m.
We use a risk aversion parameter $\gamma=\num{1e-2}$.

The optimal policy has the form
\[
u_t = K_tW + k_t,
\]
where $K_t\in\reals^{3\times 1}$ and $k_t\in\reals^{3}$.
Therefore the consumption amount during year $t$ is equal to
\BEAS
C_t &=& W_t-\ones^T u_t\\
&=& W_t - \ones^T (K_tW_t + k_t)\\
&=& (1-\ones^T K_t) W_t - \ones^T k_t.
\EEAS
We can rewrite this in the more interpretable form
\[
u_t= f_t (W_t-2000) + g_t.
\]
Here $(W_t-2000)$ is the deficit or excess wealth the investor has as compared to the optimal bequest amount, $f_t$ is unit-less, and $g_t$ is in units thousands of dollars.

\begin{figure}
    \begin{minipage}{.5\textwidth}
        \centering
        \includegraphics[width=\linewidth]{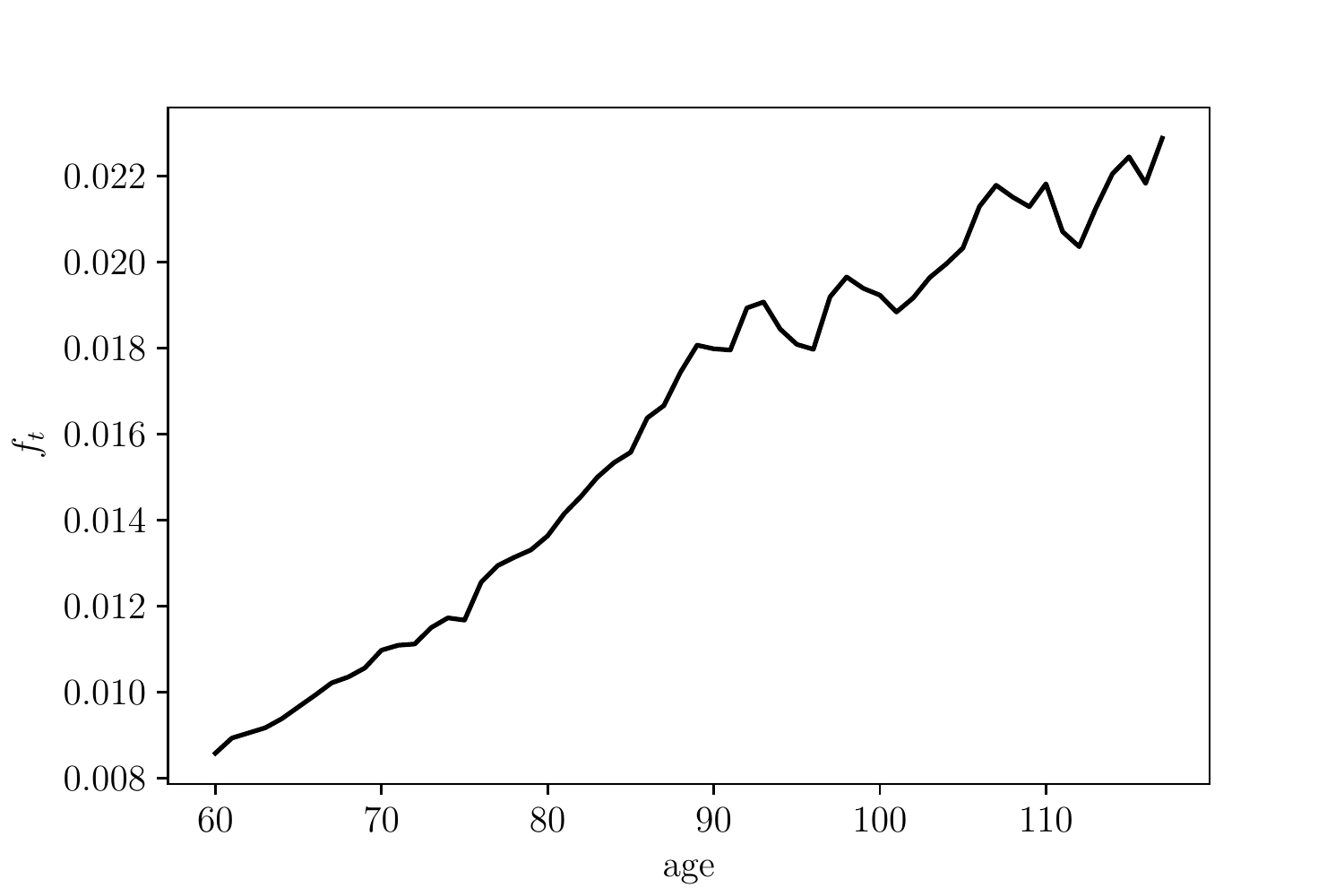}
        \caption{$f_t$ versus age.}
        \label{fig:f}
    \end{minipage}%
    \begin{minipage}{0.5\textwidth}
        \centering
        \includegraphics[width=\linewidth]{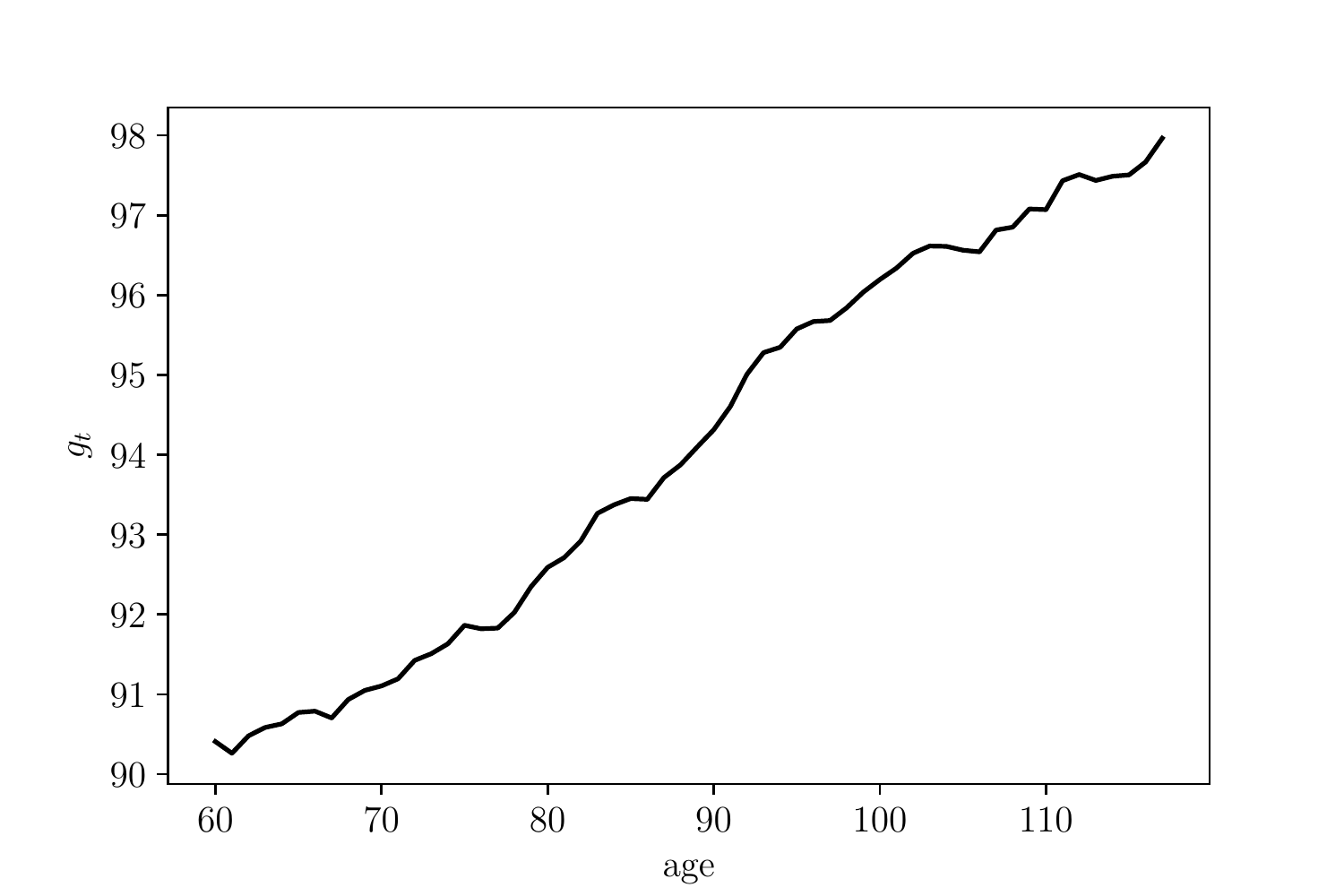}
        \caption{$g_t$ versus age.}
        \label{fig:g}
    \end{minipage}
\end{figure}

Figure~\ref{fig:f} shows $f_t$ versus $t$ and Figure~\ref{fig:g} shows $g_t$ versus $t$, starting from age $60$.

\begin{figure}
    \centering
    \begin{minipage}{.5\textwidth}
        \centering
        \includegraphics[width=\linewidth]{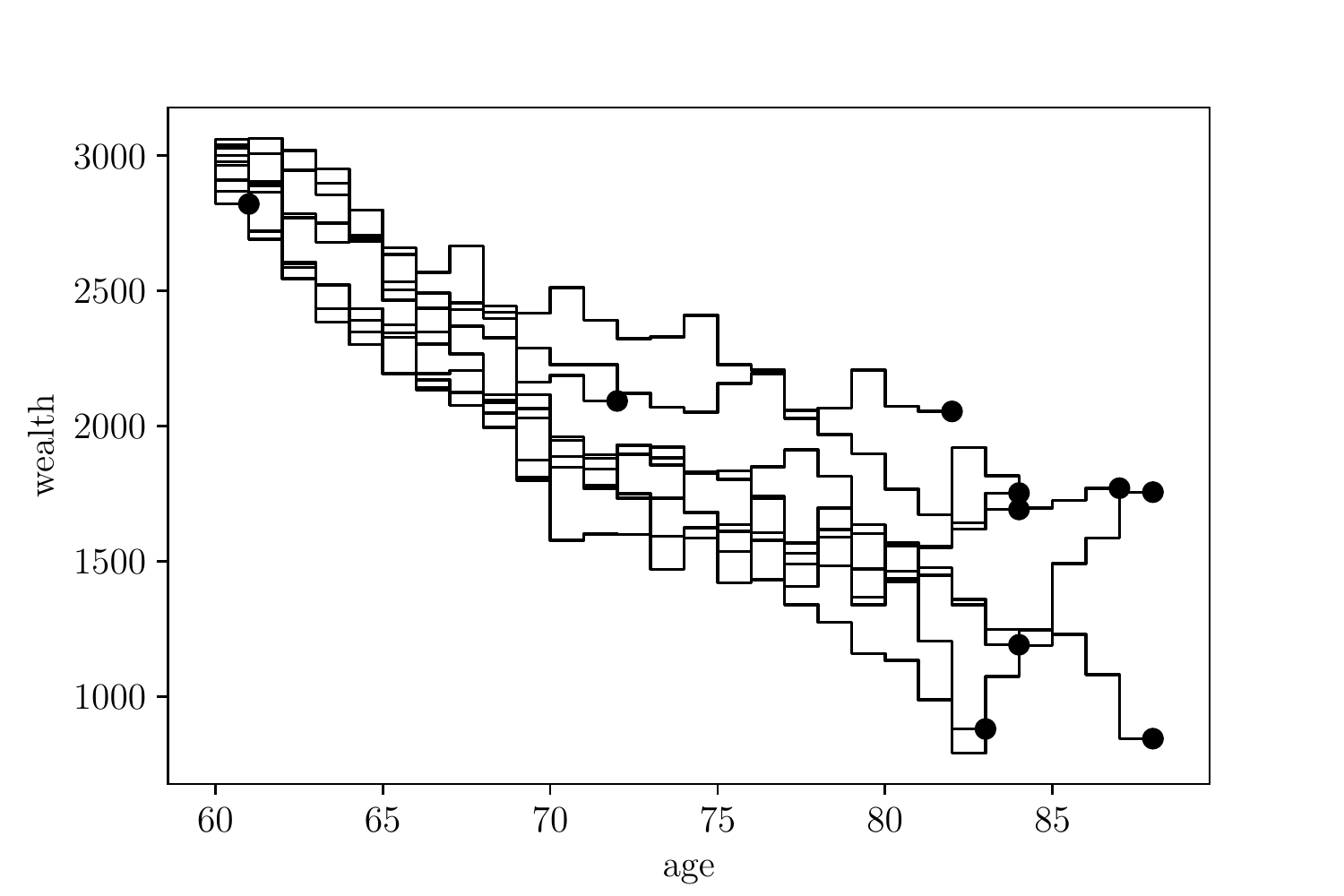}
        \caption{Wealth versus age.}
        \label{fig:wealth}
    \end{minipage}%
    \begin{minipage}{0.5\textwidth}
        \centering
        \includegraphics[width=\linewidth]{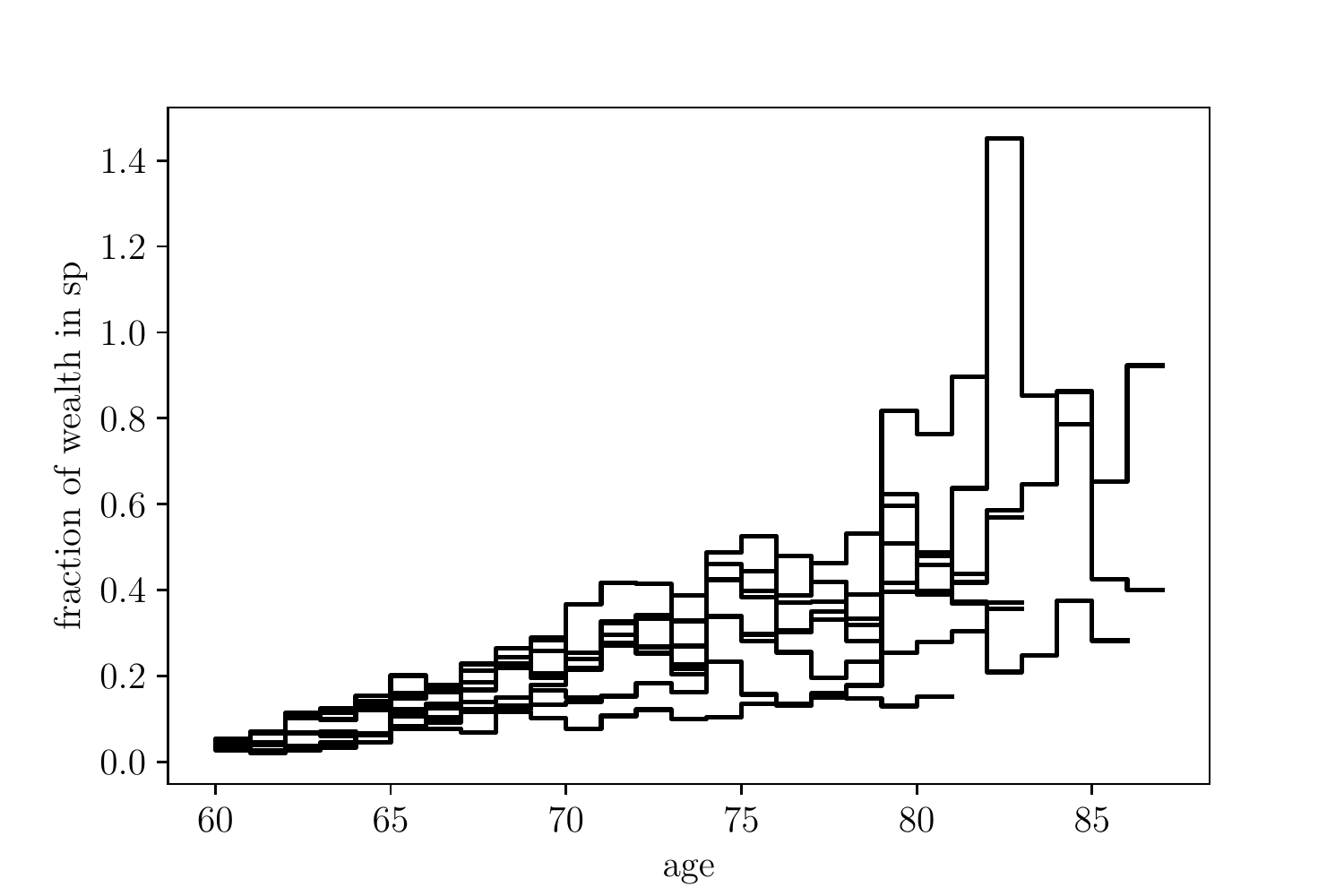}
        \caption{Fraction of wealth invested in the S\&P 500 versus age starting at age $60$ with $3$m dollars.}
        \label{fig:risky}
    \end{minipage}
\end{figure}

We simulated (the remainder of) an investor's life, starting at age $60$ with \$$3$m, using the optimal policy found by our algorithm, over $500$ random seeds.
Figure~\ref{fig:wealth} shows wealth versus age over the first $10$ random seeds.
Figure~\ref{fig:risky} shows the fraction of the investor's wealth invested in the S\&P 500 over age, over the same $10$ random seeds.
It appears that the investor seeks a riskier allocation as they age, which is counterintuitive, since one would expect the opposite.

\begin{figure}
\begin{minipage}{0.5\textwidth}
    \centering
    \includegraphics[width=\linewidth]{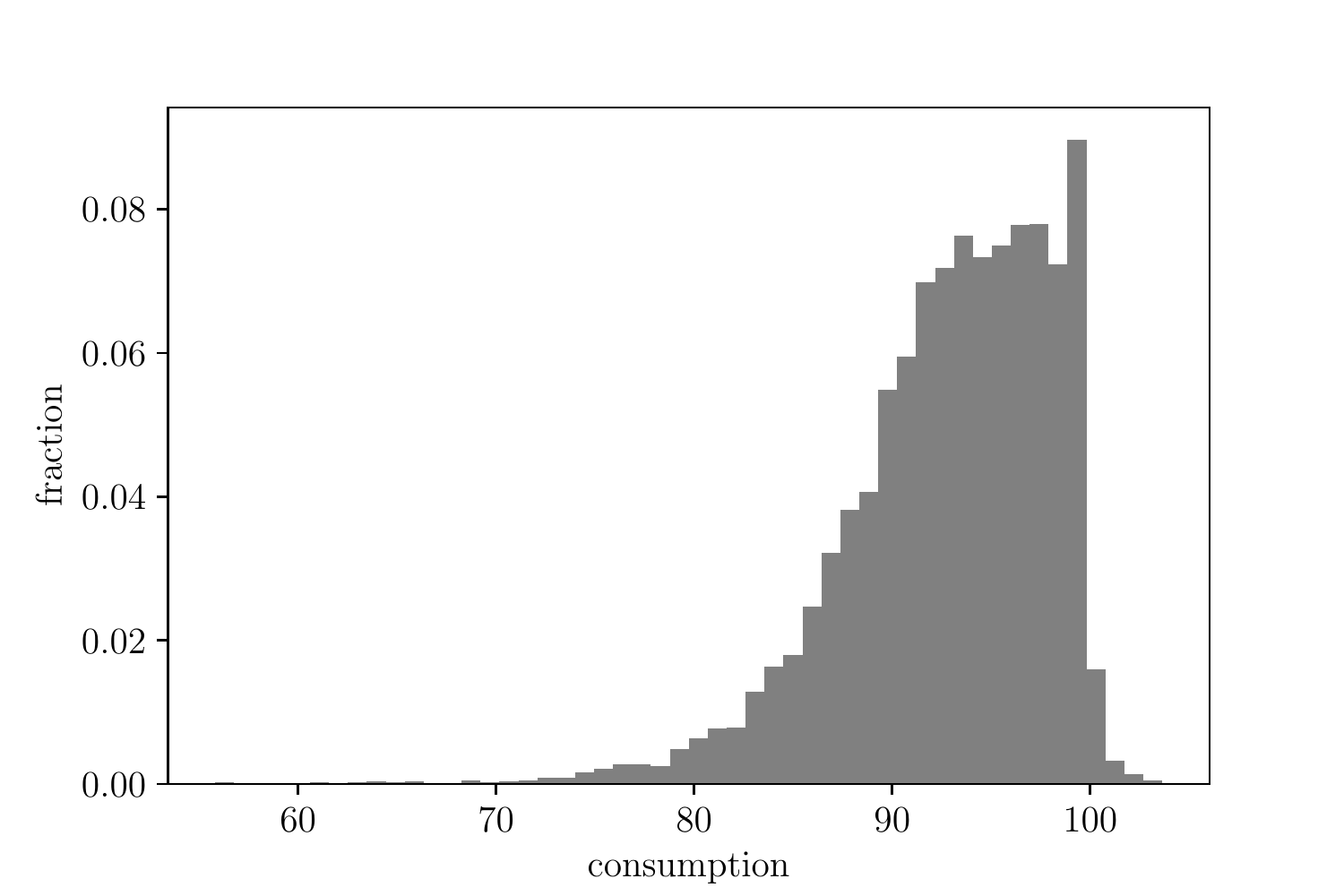}
    \caption{Histogram of consumption amounts.}
    \label{fig:consumption_hist}
\end{minipage}
\begin{minipage}{0.5\textwidth}
    \centering
    \includegraphics[width=\linewidth]{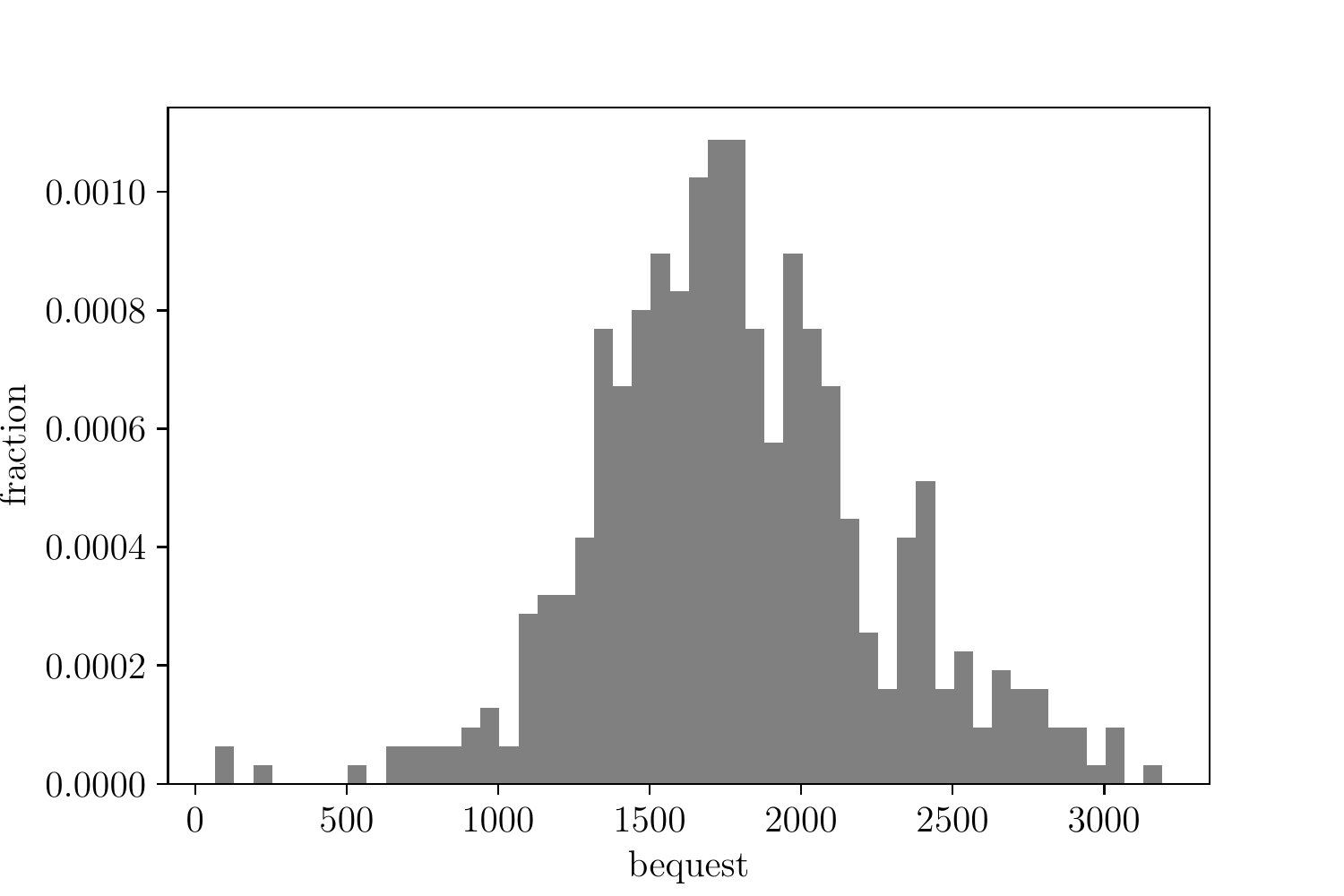}
    \caption{Histogram of bequest amounts.}
    \label{fig:bequest_hist}
\end{minipage}
\end{figure}

Figure~\ref{fig:consumption_hist} shows a histogram of consumption amounts, over all $500$ random seeds.
The investor rarely consumes over $100$k; if they do, this is likely because they do not want to bequest \emph{too much} (this is a limitation of quadratic utility functions).
Figure~\ref{fig:bequest_hist} shows a histogram of bequest amounts, over all $500$ random seeds.
The median bequest amount was \$$1.747$m.

\section{Conclusion}

It has been about sixty years since the invention of LQR.
Since its invention, everyone has known that LQR problems can be solved exactly.
Many extensions, some which maintain tractability, and some that do not, have been proposed over the years.
In this paper, we have collected the tractable extensions, unified them as a single general class of problems,
and proven the form of the solution.
There is no clean expression for the optimal policy, however computing it can be reduced to iteratively performing several
simple linear algebraic operations on the coefficients of extended quadratic functions.
We have also developed an implementation which exactly solves (modulo how expectation is performed) these problems
using these operations defined on extended quadratics, given access only to a sampling oracle.
We demonstrate the usefulness of such an approach via many applications, some of which, to the best of our knowledge, have
not appeared yet in the literature. 

\clearpage
\section*{Acknowledgements}
Shane Barratt is supported by the National Science Foundation Graduate Research Fellowship under Grant No. DGE-1656518.

\bibliography{gen_lqr}

\appendix
\newpage
\section{Exact expectation given first and second moments}
\label{sec:exact}
\paragraph{Expectation of affine pre-composition.}
Consider the (non-extended) quadratic function
\[
g(x) = \frac{1}{2}\begin{bmatrix}x \\ 1\end{bmatrix}^T\begin{bmatrix} P & q \\ q^T & r\end{bmatrix}\begin{bmatrix}x \\ 1\end{bmatrix},
\]
where $x\in\reals^n$, with deterministic coefficients.
Suppose we are given the random matrix $K\in\reals^{n \times m}$ and random vector $k\in\reals^m$ with a joint distribution.
Then
\[
h(z) = \Expect g(Kz+k)
\]
is a (non-extended) quadratic function of $z$.
We can express $h(z)$ as
\[
h(z) =
\frac{1}{2}
\begin{bmatrix}z \\ 1\end{bmatrix}^T
\Expect
\left[
A^T
B
A
\right]
\begin{bmatrix}z \\ 1\end{bmatrix}
\]
where $A=\begin{bmatrix}
        K & k \\
        0 & 1
    \end{bmatrix}$ and $B=\begin{bmatrix}
        P &q \\
        q^T & r
    \end{bmatrix}$.
Then we can write the $ij$th entry of the above expectation as
\[
\Expect[A^TBA]_{ij} = \sum_l \sum_k B_{lk}\Expect[A_{li}A_{kj}].
\]
So if we know the first and second moments of $(K,k)$, which in turn give us the second moments of $A$, we can calculate the resulting (non-extended) quadratic function analytically.

\paragraph{Expectation of random quadratic.}
Consider the (non-extended) quadratic function
\[
g(x) = \frac{1}{2}\begin{bmatrix}x \\ 1\end{bmatrix}^T\begin{bmatrix} P & q \\ q^T & r\end{bmatrix}\begin{bmatrix}x \\ 1\end{bmatrix},
\]
where $x\in\reals^n$, with random coefficients.
Then
\[
\Expect g(x) = \frac{1}{2}\begin{bmatrix}x \\ 1\end{bmatrix}^T\Expect\begin{bmatrix} P & q \\ q^T & r\end{bmatrix}\begin{bmatrix}x \\ 1\end{bmatrix}
\]
can trivially be computed with knowledge of the first moment of the coefficients.

\section{Lognormal distribution}
\label{sec:lognormal}

Suppose $x\in\reals^n$ is a random variables that follows a lognormal distribution, \ie,
\[x=\exp(z),\]
where
$z\sim\mathcal{N}(\mu,\Sigma)$.

The mean of $x$ is
\[
\Expect [x] = e^{\mu + \frac{1}{2}\diag(\Sigma)}
\]
and the covariance of $x$ is
\[
\text{Cov} [x] = \diag(\Expect[x])(e^{\Sigma}-\ones\ones^T)\diag(\Expect[x])
\]
where $(e^\Sigma)_{ij} = e^{\Sigma_{ij}}$.

\section{Code example}

The following python code is an example of the code one would write to solve a random instance of the infinite-horizon LQR problem:
\begin{verbatim}
As = np.random.randn(1,1,n,n)
Bs = np.random.randn(1,1,n,m)
cs = np.zeros((1,1,n))
gs = [ExtendedQuadratic(np.eye(n+m),np.zeros(n+m),0)]
Pi = np.eye(1)
def sample(N):
    A = np.zeros((N,K,n,n)); A[:] = As
    B = np.zeros((N,K,n,m)); B[:] = Bs
    c = np.zeros((N,K,n)); c[:] = cs
    g = [gs for _ in range(N)]
    return A,B,c,g,Pi
Vs, Qs, policies = dp_infinite(sample, num_iterations=50, N=1, gamma=1)
\end{verbatim}

\end{document}

%% file: defs.tex
\newcommand{\BEAS}{\begin{eqnarray*}}
\newcommand{\EEAS}{\end{eqnarray*}}
\newcommand{\BEQ}{\begin{equation}}
\newcommand{\EEQ}{\end{equation}}
\newcommand{\BIT}{\begin{itemize}}
\newcommand{\EIT}{\end{itemize}}

\newcommand{\eg}{{\it e.g.}}
\newcommand{\ie}{{\it i.e.}}

\newcommand{\ones}{\mathbf 1}
\newcommand{\reals}{{\mbox{\bf R}}}

\newcommand{\symm}{{\mbox{\bf S}}}  

\newcommand{\nullspace}{{\mathcal N}}
\newcommand{\range}{{\mathcal R}}
\newcommand{\Rank}{\mathop{\bf rank}}

\newcommand{\diag}{\mathop{\bf diag}}

\newcommand{\Expect}{\mathop{\bf E{}}}

\newcommand{\argmin}{\mathop{\rm argmin}}

\newcommand{\etal}{{\it et al.}}

\newcounter{algorithmctr}
\renewcommand{\thealgorithmctr}{\arabic{algorithmctr}}
\newenvironment{algdesc}%
   {\refstepcounter{algorithmctr}\begin{list}{}{%
       \setlength{\rightmargin}{0\linewidth}%
       \setlength{\leftmargin}{.05\linewidth}}%
       \rmfamily\small
       \item[]{\setlength{\parskip}{0ex}\hrulefill\par%
        \nopagebreak{\bfseries\textsf{Algorithm \thealgorithmctr~}}}}%
   {{\setlength{\parskip}{-1ex}\nopagebreak\par\hrulefill} \end{list}}